\documentclass[a4j,12pt]{article} 
\usepackage{amsmath,amssymb}
\usepackage{mathdots}   
\newcommand{\bfd}{{\bf d}}

\allowdisplaybreaks[1]   

\newcommand{\Section}[1]{%
\renewcommand{\thesection}{\S\arabic{section}}
\section{#1}
\renewcommand{\thesection}{\arabic{section}}
\setcounter{equation}{0}}

\newcommand{\qed}{\hfill \hbox{\rule[-2pt]{4pt}{7pt}}}
\newcommand{\proof}{{\hspace*{0.4cm} {\it Proof}.\ \enskip}}

\newtheorem{lem}{Lemma}[section]
\newtheorem{thm}[lem]{Theorem}  
\newtheorem{prop}[lem]{Proposition}
\newtheorem{cor}[lem]{Corollary}
\newtheorem{rem}[lem]{Remark}

\newcommand{\Z}{{\Bbb Z}}
\newcommand{\R}{{\Bbb R}}
\newcommand{\C}{{\Bbb C}}

\newcommand{\N}{{\Bbb N}}

\newcommand{\calH}{{\cal H}}

\newcommand{\calO}{{\cal O}}

\newcommand{\calR}{{\cal R}}
\newcommand{\calS}{{\cal S}}
\newcommand{\calU}{{\cal U}}

\newcommand{\pair}[2]{\left\langle {#1}, \, {#2}\right\rangle}

\newcommand{\set}[2]{\left\{\left.#1\vphantom{#2}\:\right\vert\:#2\right\}}
\newcommand{\wt}{\widetilde}
\newcommand{\what}{\widehat}
\newcommand{\xx}{{\frak X}}
\newcommand{\fra}{{\frak a}}

\newcommand{\frp}{{\frak p}}

\newcommand{\lam}{{\lambda}}
\newcommand{\Lam}{{\Lambda}}

\newcommand{\alp}{{\alpha}}  
\newcommand{\vphi}{{\varphi}} 
\newcommand{\eps}{{\epsilon}} 

\newcommand{\slit}{\vspace{5mm}}
\newcommand{\mslit}{\vspace{3mm}}
\newcommand{\real}{{\rm Re}}

\newcommand{{\any}}{{}^\forall}
\newcommand{{\is}}{{}^\exists}
\newcommand{{\st}}{\; {\rm s.t.}\; }

\newcommand{\ol}[1]{\overline{#1}}

\newcommand{\hec}{{\calH(G,K)}}

\newcommand{\SKX}{{\calS(K \backslash X)}}

\newcommand{\abs}[1]{\left\vert{#1}\right\vert}  

\newcommand{\dint}[1]{\displaystyle{\int_{#1}}}

\newcommand{\dprod}[1]{\displaystyle{\prod_{#1}}}

\newcommand{\dsqcup}[1]{\displaystyle{\bigsqcup_{#1}}}

\newcommand{\twomatrix}[4]{\begin{pmatrix}
                           {#1} & {#2}\\
                           {#3} & {#4}
                          \end{pmatrix}}

\newcommand{\twomatrixplus}[4]{\left(\begin{array}{c|c}
                           {#1}  & {#2}\\
                           \hline
                           {#3} &  {#4}
                          \end{array}\right)}

\newcommand{\mapright}[1]{\displaystyle{
   \smash{\mathop{\hbox to 1cm{\rightarrowfill}}^{#1}}}}

\newcommand{\gyaddots}{%
\setlength{\unitlength}{1mm}
\begin{picture}(5,3.5)(-2,-1.5)
\put(0,0){$\cdot$}
\put(2,1.5){$\cdot$}
\put(-2,-1.5){$\cdot$}
\end{picture}}
\newcommand{\gyakuddots}{\smash{\lower0.3ex\hbox{\gyaddots}}}

\begin{document}
\title{Harmonic analysis on the space of \\$p$-adic unitary hermitian matrices,  \\
including dyadic case\footnote{ 
This research is partially supported by Grant-in-Aid for Scientific Research (C): 24540031
  }  }

\author{Yumiko Hironaka} 
\date{}
\maketitle

\bigskip
\begin{minipage}{15cm}
Abstract: We are interested in the harmonic analysis on $p$-adic homogeneous spaces based on spherical functions. In the present paper, we investigate the space $X$ of unitary hermitian matrices of size $m$ over a $\frp$-adic field $k$ and give unified description including dyadic case, which is a continuation of our previous papers on non-dyadic case. The space becomes complicated when $e = v_\pi(2) > 0$. 
First we introduce a typical spherical function $\omega(x;z)$ on $X$, and study their functional equations, which depend on $m$ and $e$, we give an explicit formula for $\omega(x;z)$, where Hall-Littlewood polynomials of type $C_n$ appear as a main term with different specialization according as $m = 2n$ or $2n+1$, but independent of $e$. 
By spherical transform, we show the Schwartz space $\SKX$ is a free Hecke algebra $\hec$-module of rank $2^n$, and give parametrization of all the spherical functions on $X$ and the explicit Plancherel formula on $\SKX$. The Plancherel measure does not depend on $e$, but the normalization of $G$-invariant measure on $X$ depends. 



\end{minipage}

\setcounter{section}{-1}

\Section{Introduction}
We have been interested in harmonic analysis on $p$-adic homogeneous spaces based on spherical functions. We have considered the space of $p$-adic unitary hermitian matrices of even size in \cite{HK} and \cite{HK2} mainly for odd residual case, and in the present paper we will give the unified description including even residual case. All the results for odd residual case have been proved in \cite{HK} (resp. \cite{HK2}) for even (resp. odd) size
matrices. 
When the matrix size is even, the space 
has a natural close relation to the theory of automorphic functions and classical theory of sesquilinear forms (e.g. \cite{Oda}, \cite{Arakawa}), where we didn't need to distinguish the dyadic case. 

We fix an unramified quadratic extension $k'$ of $\frp$-adic field $k$, and consider hermitian and unitary matrices with respect to $k'/k$, and for $a \in M_{mn}(k')$ we denote by $a^*\in M_{nm}(k')$ its conjugate transpose.  
Let $\pi$ be a prime element of $k$ and $q$ the cardinality of the residue class field $\calO_k/(\pi)$. 
Denote by $j_m \in GL_m(k)$ the matrix whose all anti-diagonal entries are $1$ and others are $0$, where and henceforth $m$ is an integer such that $m \geq 2$. Set
\begin{eqnarray*}
&&
G = U(j_m) = \set{g\in GL_m(k')}{g^*j_mg = j_m}, \\ 
&&
X = \set{x \in X}{x^* = x, \; \Phi_{xj_m}(t) = \Phi_{j_m}(t)},\\
&&
g \cdot x = gxg^*, \quad (g \in G, \; x \in X),
\end{eqnarray*}
where $\Phi_y(t)$ is the characteristic polynomial of matrix $y$. We note that $X$ is a single $G(\ol{k})$-orbit over the algebraic closure $\ol{k}$ of $k$. 
We fix $K = G \cap GL_m(\calO_{k'})$, which is a maximal compact open subgroup of $G$ satisfying the Iwasawa decomposition $G = KB = BK$ with Borel subgroup $B$ of $G$ consisting of all the upper triangular matrices in $G$. 

For dyadic case, i.e. $v_\pi(2) >0$, there are $K$-orbits in $X$ which have no diagonal element (cf. Theorem~1 below); while for the space of unramified hermitian matrices in $GL_m(k')$, each $GL_m(\calO_{k'})$-orbit has a diagonal representative (\cite{Jac}).   
It is known in general that the spherical functions on various $p$-adic groups $\Gamma$ can be written in terms of the specialization of Hall-Littlewood polynomials of the corresponding root structure of $\Gamma$ (cf.\ \cite[\S 10]{Mac}, also \cite[Theorem 4.4]{Car}).  For the present space $X$, the main term of spherical functions can be written by using Hall-Litllewood polynomials of type $C_n$ with different specialization according to the parity of $m$, independent of the residual characteristic (cf.~Theorem~3 below). By using spherical functions we study the Schwartz space $\SKX$, and show its $\hec$-module structure, parametrization of all the spherical functions on $X$ and Plancherel formula and Inversion formula on $\SKX$(cf.~Theorem~4 below).  

\bigskip
We will explain the results in some more details. Set
\begin{eqnarray} \label{eq:m-n}
n = \left[\frac{m}{2}\right], \qquad e = v_\pi(2) \, (\geq 0)
\end{eqnarray}
and denote the corresponding groups $G, \; B, \; K$ and space $X$ with subscript and superscript  if necessary, as  $G_n^{(ev)},  \; X_n^{(ev)}$ for $m = 2n$, and  $G_n^{(od)},  \; X_n^{(od)}$ for $m = 2n+1$, etc.
According to the parity of $m$, $G$ has the root structure of type $C_n$ for even $m$ or type $BC_n$ for  odd $m$. We fix a unit $\eps \in k$ for which $\{1, \frac{1+\sqrt{\eps}}{2}\}$ forms an $\calO_k$-basis for $\calO_{k'}$, where $v_\pi(1-\eps)= 2e$ (cf. \cite[\S 64]{Ome}).
Set
\begin{eqnarray}
&& \label{eq:wt-Lam}
\wt{\Lam_n^+} = \set{\lam \in \Z^n}{\lam_1 \geq \lam_2 \geq \cdots \geq \lam_n \geq -e}, \\
&&
\Lam_n^+ = \set{\lam \in \wt{\Lam_n^+}}{\lam_n \geq 0}, \nonumber
\end{eqnarray}
where $\Lam_n^+ = \wt{\Lam_n^+}$ if $e = 0$. 
For $\lam \in \wt{\Lam_n^+}$ such that $\lam_r \geq 0 > \lam_{r+1}$, define
\begin{eqnarray*}
&&
x_\lam^{(ev)} = Diag(\pi^{\lam_1}, \ldots, \pi^{\lam_r}, y_\lam^{(ev)}, \pi^{-\lam_r}, \ldots, \pi^{-\lam_1}) \in X_n^{(ev)}, \\
&&
\quad y_\lam^{(ev)} = \left\{\begin{array}{ll}
\emptyset & \mbox{if  $r = n$}\\
\begin{pmatrix}
    \pi^{\lam_{r+1}}(1-\eps) &&&&& -\sqrt{\eps} \\
    & \ddots &&& \iddots & \\
    &&\pi^{\lam_n}(1-\eps)& -\sqrt{\eps}&& \\
&&\sqrt{\eps}&\pi^{-\lam_n} && \\
    & \iddots &&& \ddots & \\
    \sqrt{\eps} &&&&& \pi^{-\lam_{r+1}}
  \end{pmatrix} & \mbox{if $r<n$},
\end{array} \right.  \nonumber \\
&&
x_\lam^{(od)} = Diag(\pi^{\lam_1}, \ldots, \pi^{\lam_r}, y_\lam^{(od)}, \pi^{-\lam_r}, \ldots, \pi^{-\lam_1}) \in X_n^{(od)},\\
&&
\quad y_\lam^{(od)} = \left\{\begin{array}{ll}
1 & \mbox{if  $r = n$}\\
\begin{pmatrix}
    \pi^{\lam_{r+1}}(1-\eps) &&&&&& -\sqrt{\eps} \\
    & \ddots &&&& \iddots & \\
    &&\pi^{\lam_n}(1-\eps)&& -\sqrt{\eps}&& \\
&& & 1 &&&\\    
&&\sqrt{\eps}&&\pi^{-\lam_n} && \\
    & \iddots &&&& \ddots & \\
    \sqrt{\eps} &&&&&& \pi^{-\lam_{r+1}}
  \end{pmatrix} & \mbox{if $r<n$},
\end{array} \right. \nonumber
\end{eqnarray*}
and understand $x_\lam = y_\lam$ if $r = 0$, where and henceforth we write simply as $x_\lam$ or $y_\lam$, if there is no danger of confusion. Here and henceforth empty entries in matrices should be understood as $0$. %

\slit
\noindent
{\bf Theorem~1} (1)  {\it The map $\wt{\Lam_n^+} \longrightarrow K \backslash X, \; \lam \longmapsto K\cdot x_\lam$ is surjective. Further, it is bijective if $m$ is even, $m = 3$, or $e =1$. }

\noindent
(2) {\it There are precisely two $G$-orbits in $X$ represented by $x_\lam$ with $\lam = \bf0$ and $(1,0,\ldots,0)$. }

\slit
It is known that  each $x_\lam, \; \lam \in \Lam_n^+$,  gives a different $K$-orbit, since it gives a different $GL_m(\calO_{k'})$-orbit in the space of hermitian matrices in $GL_m(k')$, where $m = 2n$ or $2n+1$. Hence, when $e = 0$, it is enough to show that any $K$-orbit has a representative of the shape $x_\lam, \; \lam \in \Lam_n^+$, which has been done in \cite{HK} for even $m$ and \cite{HK2} for odd $m$.  For $e >0$, we have to prove also the non-redundancy within the above representatives, which will be done as a corollary of the explicit formula of spherical functions on $X$  for $n \geq 2$, i.e., $m \geq 4$ (cf. Theorem~3).

\mslit
A spherical function on $X$ is a $K$-invariant function on $X$ which is a common eigenfunction with respect to the convolutive action of the Hecke algebra $\hec$, and a typical one is constructed by Poisson transform from relative invariants of a parabolic subgroup.  We take the Borel subgroup $B$ consisting of upper triangular matrices in $G$.  
For $x \in X$ and $s \in \C^n$, we consider the following integral
\begin{eqnarray} \label{def-sph0}
\omega(x; s) = \int_{K} \prod_{i=1}^n \abs{d_i(k\cdot x)}^{s_i} dk,
\end{eqnarray}
where $\abs{\;}$ is the absolute value on $k$ normalized by $\abs{\pi} = q^{-1}$, $q = \sharp\left(\calO_k/(\pi)\right)$, $d_i(y)$ is the determinant of the lower right $i$ by $i$ block of $y$, $1 \leq i \leq n$, and $dk$ is the normalized Haar measure on $K$. 
Then the right hand side of (\ref{def-sph0}) is absolutely convergent for $\real(s_i) \geq 0, \; 1 \leq i \leq n$, and continued to a rational function of $q^{s_1}, \ldots , q^{s_n}$, and we use the notation $\omega(x; s)$ in such sense. Since $d_i(x)$'s are  relative $B$-invariants on $X$ such that 
\begin{eqnarray*}
d_i(p \cdot x) = \psi_i(p)d_i(x), \quad \psi_i(p) = N_{k'/k}(d_i(p)) \quad (p \in B, \; x \in X,\; 1 \leq i \leq n),
\end{eqnarray*}
we see $\omega(x; s)$ is a spherical function on $X$ which satisfies
\begin{eqnarray*}
&&
f * \omega(x; s) = \lam_s(f) \omega(x; s), \quad f \in \hec\\
&&
\lam_s(f) = \int_{B}f(p)\prod_{i=1}^n\abs{\psi_i(p)}^{-s_i} \delta(p)dp,
\end{eqnarray*}
where $dp$ is the normalized left invariant measure on $B$ with modulus character $\delta$.  
The Weyl group $W$ of $G$ relative to $B$ acts on rational characters of $B$, hence on $s$.
It is convenient to introduce the new variable $z \in \C^n$ related to $s$ by
\begin{eqnarray}
&&
s_i = -z_i + z_{i+1} -1 + \frac{\pi\sqrt{-1}}{\log q}, \quad 1 \leq i \leq n-1 \nonumber \\
&&
\label{eq:chgv}
s_n = \left\{\begin{array}{ll} 
     -z_n -\frac12 & \mbox{if }m = 2n \\[2mm]
     -z_n -1+ \frac{\pi\sqrt{-1}}{2\log q} & \mbox{if } m = 2n+1. 
\end{array}\right.
\end{eqnarray}
Then $W = \langle S_n, \tau \rangle$ acts on $z$ by permutation of indices as for the elements of $S_n$ and  $\tau(z_1, \ldots, z_n)= (z_1, \ldots, z_{n-1}, -z_n)$. 
Keeping the above relation \eqref{eq:chgv}, we denote $\omega(x;z) = \omega(x;s)$ and $\lam_s = \lam_z$. Then $\lam_z$ gives the Satake isomorphism 
\begin{eqnarray} \label{eq:Satake}
\lam_z: \hec \stackrel{\sim}{\longrightarrow} \C[q^{\pm 2z_1}, \ldots, q^{\pm 2z_n}]^W (= \calR_0, \mbox{say}).
\end{eqnarray}
We will give the functional equation of $\omega(x;z)$ with respect to $W$. 
We set 
\begin{eqnarray*}
&&
\Sigma^+ = \Sigma_s^+ \sqcup \Sigma_\ell^+,\\
&&
\Sigma_s^+ = \set{e_i + e_j, \; e_i -e_j}{1 \leq i < j \leq n}, \quad 
\Sigma_\ell^+ = \set{2e_i}{1 \leq i  \leq n}, 
\end{eqnarray*}
where $e_i \in \Z^n$ is the $i$-th unit vector, and define a pairing
\begin{eqnarray*}
\pair{\;}{} : \Z^n \times \C^n \longrightarrow \C, \; \pair{\alp}{z} = \sum_{i=1}^n \alp_iz_i.
\end{eqnarray*}

\bigskip
\noindent
{\bf Theorem~2} \textit{Assume $e \leq 1$ if $m$ is odd.}\\
{\rm (1)} \textit{For any $\sigma \in W$, one has}
$$
\omega(x; z) = \Gamma_\sigma^{(e)}(z)  \cdot \omega(x; \sigma(z)).
$$
where
$$
\Gamma_\sigma^{(e)}(z) = \prod_{\alp}\, \gamma_\alp^{(e)}(z), \quad
\gamma_\alp^{(e)}(z) = \left\{\begin{array}{ll}
\frac{1-q^{-1+\pair{\alp}{z}}}{q^{\pair{\alp}{z}}-q^{-1}} & \alp \in \Sigma_s^+ \\[2mm]
q^{e\pair{\alp}{z}} & \alp \in \Sigma_\ell^+, m=2n \\[2mm]
q^{e\pair{\alp}{z}}\frac{1-q^{-1+\pair{\alp}{z}}}{q^{\pair{\alp}{z}}-q^{-1}} & \alp \in \Sigma_\ell^+, m=2n+1,
\end{array} \right.
$$
\textit{and  $\alp$ runs over the set $\set{\alp \in \Sigma^+}{-\sigma(\alp) \in \Sigma^+}$.}

\medskip
\noindent
{\rm (2)} \textit{The function $q^{-\pair{e}{z}}G(z) \cdot \omega(x;z)$ is holomorphic and $W$-invariant, hence belongs to $\C[q^{\pm z_1}, \ldots, q^{\pm z_n}]^W$. Here
$$
\pair{e}{z} = e(z_1+\cdots +z_n), \quad G(z) = \prod_{\alp}\, \frac{1+q^{\pair{\alp}{z}}}{1-q^{-1+\pair{\alp}{z}}}, 
$$
and $\alp$ runs over the set $\Sigma_s^+$ for $m = 2n$ and $\Sigma^+$ for $m = 2n+1$.}

\bigskip
As for the explicit formula of $\omega(x; s)$ it suffices to give for $x_\lam$ by Theorem 1 (1). 

\bigskip
\noindent
{\bf Theorem~3} {\rm (Explicit formula)} {\it Assume $e \leq 1$ if $m$ is odd. For each $\lam \in \wt{\Lam_n^+}$, one has }
\begin{eqnarray*} 
\omega(x_\lam; z)
 = c_n\, q^{\pair{\lam}{z_0}} \cdot \frac{q^{\pair{e}{z}}}{G(z)}  \cdot Q_{\lam+e}(z; \{t\}),
\end{eqnarray*}
\textit{where $G(z)$ is given in Theorem 2,  $z_0$ is the value in $z$-variable corresponding to $s = {\bf 0}$,}
\begin{eqnarray*}
&&
c_n = \left\{\begin{array}{ll}
\dfrac{(1-q^{-2})^n}{w_m(-q^{-1})} & \textit{if } m = 2n\\
\dfrac{(1+q^{-1})(1-q^{-2})^n}{w_m(-q^{-1})} & \textit{if } m = 2n+1,
\end{array} \right. \quad w_m(t) = \prod_{i=1}^m (1-t^i),\\ 
&&
Q_\mu(z; \{t\}) = \sum_{\sigma \in W}\, \sigma\left(q^{-\pair{\mu}{z}} c(z, \{t\}) \right),\quad
c(z,\{t\}) = \prod_{\alp \in \Sigma^+}\, \frac{1 -t_\alp q^{\pair{\alp}{z}}}{1 - q^{\pair{\alp}{z}}},\\
&& 
\{t\} = \{t_\alp\} \quad \mbox{with}\; \;  t_\alp = \left\{ \begin{array}{ll} 
-q^{-1} & \textit{if } \alp \in \Sigma_s^+ \\[2mm] 
q^{-1} & \textit{if } \alp \in \Sigma_\ell^+, \; m = 2n\\[2mm]
-q^{-2} & \textit{if } \alp \in \Sigma_\ell^+, \; m = 2n+1. \end{array}\right.
\end{eqnarray*}

\bigskip
We see the main part $Q_{\lam+e}(z; \{t\})$ of $\omega(x_\lam;z)$ belongs to $\calR = \C[q^{\pm z_1}, \ldots, q^{\pm z_n}]^W$ by Theorem 2.  It is known that $Q_\mu(z;\{t\})$ is a constant multiple of Hall-Littlewood polynomial $P_\mu(z; \{t\})$ and the set $\set{P_\mu(z; \{t\})}{\mu \in \Lam_n^+}$ forms a $\C$-basis for $\calR$ (for more details, see Remark~\ref{Rem 3-1}). Hence we see the non-redundancy for the representatives in Theorem~1-(1). 
The influence of the residual characteristic of the base field $k$ in the explicit formula of $\omega(x_\lam; z)$ appears as shifting $\lam+e$ in $Q_{\lam}$ or $P_{\lam}$ and the factor $q^{\pair{e}{z}}$. 

We modify the spherical function by using the value at $x_{(-e)}, \; (-e) \in \wt{\Lam_n^+}$ as 
\begin{eqnarray*} \label{modified sph}
\Psi(x; z) = \frac{\omega(x; z)}{\omega(x_{(-e)}; z)} \in \C[q^{\pm z_1}, \ldots, q^{\pm z_n}]^W (= \calR), 
\end{eqnarray*}
and define the spherical Fourier transform on the Schwartz space $\SKX$ by
\begin{eqnarray*}
&&
\what{\; } :  \SKX  \longrightarrow  \calR, \;  \vphi  \longmapsto  \what{\vphi}(z) = \int_{X} \vphi(x) \Psi(x;z) dx
\end{eqnarray*}
where $dx$ is a $G$-invariant measure on $X$. Then it satisfies
\begin{eqnarray*}
(f*\vphi)^{\what{\;}} = \lam_z(f) \what{\vphi}, \quad f \in \hec, \; \vphi \in \SKX,
\end{eqnarray*}
where $\lam_z$ is the Satake isomorphism given in \eqref{eq:Satake}.

\bigskip
\noindent
{\bf Theorem~4} 
{\it Assume $e \leq 1$ if $m$ is odd.\\
{\rm (1)} The above spherical Fourier transform is an $\hec$-module isomorphism and $\SKX$ becomes a free $\hec$-module of rank $2^n$.

\noindent
{\rm (2)} All the spherical functions on $X$ are parametrized by $z \in \left(\C\big{/}\frac{2\pi\sqrt{-1}}{\log q} \right)^n \big{/} W $ through $\lam_z$, and the set $\set{\Psi(x; z + u)}{u \in \{0, \frac{\pi\sqrt{-1}}{\log q} \}^n }$ forms a $\C$-basis of spherical functions corresponding to $z$.

\noindent
{\rm (3)} {\rm (Plancherel formula)} Set a measure $d\mu(z)$ on $\fra^* = \left\{\sqrt{-1}\left(\R \big{/} \frac{2\pi}{\log q} \Z \right) \right\}^n$ by
\begin{eqnarray*}
d\mu(z) = \frac{1}{2^n n!} \cdot \frac{w_n(-q^{-1}) w_{m'}(-q^{-1})}{(1+q^{-1})^{m'}} \cdot \frac{1}{\abs{c(z,\{t\})}^2}\, dz, \quad m' = \left[\frac{m+1}{2}\right],
\end{eqnarray*}
where $dz$ is the Haar measure on $\fra^*$. By the normalized $G$-invariant measure $dx$ on $X$ {\rm (}explicitly given
 in Lemma~\ref{lem:v(Kx)}{\rm )}, one has
\begin{eqnarray*}
\int_{X} \vphi(x)\ol{\psi(x)} dx = \int_{\fra^*}\what{\vphi}(z) \ol{\what{\psi}(z)} d\mu(z) \quad (\vphi, \psi \in \SKX).
\end{eqnarray*}

\noindent
{\rm (4)} {\rm (Inversion formula)} For any $\vphi \in \SKX$, one has%
\begin{eqnarray*}
\vphi(x)  = \int_{\fra^*} \what{\vphi}(z) \Psi(x; z) d\mu(z), \quad x \in X.
\end{eqnarray*}
}

\bigskip
\noindent
The spherical function $\Psi(x;z)$ and the $G$-invariant measure $dx$ on $X$ depend on $m$ and $e = v_\pi(2)$, while the Plancherel measure $d\mu(z)$ depends only on $m$. 
A key point to establish Theorems 2, 3 and 4 for $e>0$ is the functional equation of $\omega(x;z)$ for $n = 1$, i.e., $m = 2$ and $3$ (cf. Proposition~\ref{prop: even n=1} and Proposition~\ref{prop: odd n=1} ). If \eqref{eq:fun-eq m=3} in Proposition~\ref{prop: odd n=1} is true for general $e > 0$, then we may erase the assumption $e \leq 1$ for odd $m$ in Theorems~1, 2, 3 and 4 (cf. Remark~\ref{Rem:2-4}).


\vspace{2cm}
\Section{The space $X$}   
{\bf 1.1.}
Let $k'$ be an unramified quadratic extension of a $\frp$-adic field $k$ and consider hermitian and unitary matrices with respect to $k'/k$, and for $a \in M_{mn}(k')$ we denote by $a^*\in M_{nm}(k')$ its conjugate transpose. Let $\pi$ be a prime element of $k$ and $q$ the cardinality of the residue class field $\calO_k/(\pi)$, and we normalize the absolute value on $k$ by $\abs{\pi}=q^{-1}$ and denote by $v_\pi(\; )$ the additive valuation on $k$. We set $e = v_\pi(2)$. We fix a unit $\eps \in \calO_k^\times$ for which $\{1, \frac{1+\sqrt{\eps}}{2} \}$ forms an $\calO_k$-basis for $\calO_{k'}$ (cf. \cite[\S 64]{Ome}). Then $v_\pi(1-\eps) = 2e$. 
We denote by $N$ the norm map $ N_{k'/k}$, and set
$$
j_m = \begin{pmatrix}  &&1\\& \iddots &\\ 1 & {}&  \end{pmatrix} \in M_{m},
$$
where and henceforth empty entries in matrices should be
understood as $0$.  

\bigskip
We consider the unitary group for $m \geq 2$
\begin{eqnarray*}
G = G(j_m) = \set{g \in GL_m(k')}{g^*j_mg = j_m}, \quad 
\end{eqnarray*}
and the spaces of hermitian matrices in $G$
\begin{eqnarray}
\wt{X} = \set{x \in G}{x^* = x}, \quad X = \set{x \in \wt{X}}{\Phi_{xj_m}(t) = \Phi_{j_m}(t)},
\end{eqnarray}
where $\Phi_y(t)$ is the characteristic polynomial of a matrix $y$, and we see $\det(x) = 1$ for any $x \in X$. The group $G$ acts on $\wt{X}$ and $X$ by 
\begin{eqnarray}
g\cdot x = gxg^* = gxj_mg^{-1}j_m, \quad g \in G, \; x \in \wt{X}.
\end{eqnarray}
This action can be extended to the algebraic closure $\ol{k}$ of $k$, and the set $\wt{X}(\ol{k})$ is decomposed into  $G(\ol{k})$-orbits as follows (cf. \cite[Appendix A]{HK}) : 
\begin{eqnarray} \label{eq:wtX-decomp}
\wt{X}(\ol{k}) = \bigsqcup_{i=0}^m\, \set{x \in \wt{X}(\ol{k})}{\Phi_{xj_m}(t) = (t-1)^i(t+1)^{m-i}}.
\end{eqnarray}
Then $X(\ol{k})$ is a single $G(\ol{k})$-orbit containing $1_m$ and corresponding to $i=\left[ \frac{m}{2}\right]$,  and $X = X(\ol{k}) \cap G$.  It is easy to see the following:
\begin{eqnarray}
&&   \label{eq:m=2}
\mbox{If $m=2$, then } \wt{X} = \{j_2\} \sqcup \{-j_2\} \sqcup X. \nonumber \\
&& \label{eq:m=3} 
\mbox{If $m=3$, then } \wt{X} = \{j_3\} \sqcup \{-j_3\} \sqcup X \sqcup (-X). \label{eq:m=2,3}
\end{eqnarray}
We fix a maximal compact subgroup $K$ of $G$ by
\begin{eqnarray*}
K = G \cap M_m(\calO_{k'}),
\end{eqnarray*}
(cf. \cite[\S 9]{Satake}), and take a Borel subgroup $B$ of $G$ consisting of all the upper triangular matrices in $G$. Then the group $G$ satisfies the Iwsawa decomposition $G = BK = KB$.  

We are interesting in $K$-orbit decomposition of $X$. To state the results we prepare some notation.
We set
\begin{eqnarray} \label{eq:m-n}
n = \left[\frac{m}{2}\right]
\end{eqnarray}
and denote the corresponding groups $G, \; B, \; K$ and the space $X$ with subscript and superscript  if necessary, as  $G_n^{(ev)},  \; X_n^{(ev)}$ for $m = 2n$, and  $G_n^{(od)},  \; X_n^{(od)}$ for $m = 2n+1$, etc.
We set
\begin{eqnarray}
&&
\wt{\Lam_n^+} = \set{\lam \in \Z^n}{\lam_1 \geq \lam_2 \geq \cdots \geq \lam_n \geq -e}, \quad (e = v_\pi(2)) \nonumber \\
&&
\Lam_n^+ = \set{\lam \in \wt{\Lam_n^+}}{\lam_n \geq 0} \left(= \wt{\Lam_n^+} \; \mbox{if } e = 0 \right),
\end{eqnarray}
and for each $\lam \in \wt{\Lam_n^+}$ such that $\lam_r \geq 0 > \lam_{r+1}$, define $x_\lam^{(ev)} \in X_n^{(ev)}$ for $m = 2n$ and $x_\lam^{(od)} \in X_n^{(od)}$ for $m = 2n+1$ as follows
\begin{eqnarray}
&&
x_\lam^{(ev)} = Diag(\pi^{\lam_1}, \ldots, \pi^{\lam_r}, y_\lam^{(ev)}, \pi^{-\lam_r}, \ldots, \pi^{-\lam_1}), \\
&&
\quad y_\lam^{(ev)} = \left\{\begin{array}{ll}
\emptyset & \mbox{if  $r = n$}\\
\begin{pmatrix}
    \pi^{\lam_{r+1}}(1-\eps) &&&&& -\sqrt{\eps} \\
    & \ddots &&& \iddots & \\
    &&\pi^{\lam_n}(1-\eps)& -\sqrt{\eps}&& \\
&&\sqrt{\eps}&\pi^{-\lam_n} && \\
    & \iddots &&& \ddots & \\
    \sqrt{\eps} &&&&& \pi^{-\lam_{r+1}}
  \end{pmatrix} & \mbox{if $r<n$},
\end{array} \right.  \nonumber \\
&&
x_\lam^{(od)} = Diag(\pi^{\lam_1}, \ldots, \pi^{\lam_r}, y_\lam^{(od)}, \pi^{-\lam_r}, \ldots, \pi^{-\lam_1}) ,\\
&&
\quad y_\lam^{(od)} = \left\{\begin{array}{ll}
1 & \mbox{if  $r = n$}\\
\begin{pmatrix}
    \pi^{\lam_{r+1}}(1-\eps) &&&&&& -\sqrt{\eps} \\
    & \ddots &&&& \iddots & \\
    &&\pi^{\lam_n}(1-\eps)&& -\sqrt{\eps}&& \\
&& & 1 &&&\\    
&&\sqrt{\eps}&&\pi^{-\lam_n} && \\
    & \iddots &&&& \ddots & \\
    \sqrt{\eps} &&&&&& \pi^{-\lam_{r+1}}
  \end{pmatrix} & \mbox{if $r<n$},
\end{array} \right. \nonumber
\end{eqnarray}
and understand $x_\lam = y_\lam$ if $r = 0$, where and henceforth we write simply as $x_\lam$ or $y_\lam$, if there is no danger of confusion. %
\begin{thm}  \label{th:Cartan}
{\rm (1)}  The map $\wt{\Lam_n^+} \longrightarrow K_n \backslash X_n, \; \lam \longmapsto K_n\cdot x_\lam$ is surjective.

\noindent
{\rm (2)} The above map is bijective if $m = 2n$, $m = 3$, or $e = 1$. 

\noindent
{\rm (3)} There are precisely two $G_n$-orbits in $X_n$ represented by $x_0 = 1_m$ and $x_1 = Diag(\pi, 1_{m-2}, \pi^{-1})$. For $\lam \in \wt{\Lam_n^+}$, $x_\lam \in G\cdot 1_m$ if and only if $\abs{\lam} = \sum_{i=1}^n \lam_i$ is even.
\end{thm}

We recall some classical results on unramified hermitian forms (cf. \cite{Jac}).  
The group $GL_m(k')$ acts on the space $\calH_m(k') = \set{x \in GL_m(k')}{x^* = x}$ by $g \cdot x =gxg^*$, and 
\begin{eqnarray}
\calH_m(k') &=&
\bigsqcup_{\mu \in \Lam_m} GL_m(\calO_{k'})\cdot \pi^\mu  =
GL_m(k')\cdot 1_m \, \sqcup\,  GL_m(k')\cdot \pi^{(1,0\cdots, 0)},   \label{eq:herm decomp}
\end{eqnarray}
where $\Lam_m = \set{\mu \in \Z^m}{\mu_1 \geq \cdots \geq \mu_m}$, $\pi^\mu = Diag(\pi^{\mu_1}, \ldots, \pi^{\mu_m})$, and $\pi^\mu \in GL_m(k') \cdot 1_m$  if and only if  $\abs{\mu} = \sum_{i=1}^m\mu_i$ is even.

\begin{rem} {\rm 
As for (1), we have shown for  even $m$ in \cite[\S 1]{HK}
and for odd $m$ with $e = 0$ in \cite[\S 1]{HK2}. 
In \S 1.2 (resp. \S 1.3), we will show the statement (1) for $m = 3$ (resp. general odd $m$).  
The non-redundancy of the representatives for $e = 0$ follows from \eqref{eq:herm decomp}.  
For $e > 0$, we see there are $K$-orbits without any diagonal element by Proposition~\ref{prop:j-type} below, and non-redundancy for $m = 2,  3$ follows from this and the value $\ell(x-j_m)$ (cf. \eqref{eq:def of ell} below). We will see the non-redundacy for general $m$ as a corollary of the explicit formula of spherical functions in \S 3 (See Remark~\ref{Rem 3-2}).  
The property (3) is independent of the residual characteristic and we may prove in a similar way as in \cite{HK} and \cite{HK2}, so we omit the proof. We note here the stabilizer of $G(\ol{k})$ at $x = 1_m$ is isomorphic to $U(1_n)(\ol{k}) \times U(1_{m'})(\ol{k}), \; m'=\left[\frac{m+1}{2}\right]$, and explicitly given as follows (\cite[(1.5)]{HK}, \cite[Proof of Theorem~1.1]{HK2}).
\begin{eqnarray}
&&
\set{\begin{pmatrix}a & b \\ jbj& jaj \end{pmatrix}\in GL_{2n}(\ol{k})}{a+bj, a-bj \in U(1_n)(\ol{k})}, \quad \mbox{or}\\
&&
\set{\begin{pmatrix} A & b & C\\d & f & dj\\jCj& jb&jAj\end{pmatrix} \in GL_{2n+1}(\ol{k})}{A-Cj \in U(1_n)(\ol{k}), \; \twomatrix{A+Cj}{\nu b}{\nu^* d}{f} \in U(1_{n+1})(\ol{k})}, \nonumber
\end{eqnarray}
where $j = j_n$ and $\nu \in \ol{k}$ such that $\nu\nu^* = 2$. Here, we may take $\nu$ within $k'$ if $e = 0$, while for $e > 0$, we understand ${}^*$ as an extended automorphism of $\ol{k}$.  
}
\end{rem}

\begin{prop} \label{prop:j-type}
If $x \in X \cap M_m(\calO_{k'})$ satisfies $x \equiv j_m \pmod{(\pi)}$, then the orbit $K\cdot x$ has no diagonal element.
\end{prop}

\proof
If $K\cdot x$ contains a diagonal element, it must contain $1_m$. On the other hand,  since any $k \in K$ fixes $j_m$, we have $1_m \equiv j_m \pmod{(\pi)}$, which is a contradiction.
\qed

\bigskip
For $a = (a_{ij}) \in M_{m}(k')$, $a \ne 0$, we set
\begin{eqnarray} \label{eq:def of ell}
-\ell(a) = \min\set{v_\pi(a_{ij})}{1 \leq i,j \leq m},
\end{eqnarray}
and say an entry of $a$ to be {\it minimal} if its $v_\pi$-value is $-\ell(a)$. 

\begin{lem}   \label{lem:min elem}
{\rm (1)} Let $a \in M_m(\calO_{k'})$ and $b \in M_m(k')$ such that $ab\ne 0$ and $ba \ne 0$. Then, one has
$$
\ell(ab) \leq \ell(b), \quad \ell(ba) \leq \ell(b), 
$$ 
and the equalities hold if $a \in GL_m(\calO_{k'})$. 

\noindent
{\rm (2)} For any $g \in G$, one has $\ell(g) \geq 0$, and the equality holds if and only if $g \in K$.  
\end{lem}

\proof
(1) Let $a = (a_{ij}) \in M_m(\calO_{k'})$ and $b =(b_{ij})\in M_m(k')$. Then,  we have
\begin{eqnarray*}
-\ell(ab) &=& \min \set{v_\pi(\sum_{k}\, a_{ik}b_{kj})}{1 \leq i, j \leq m}
\geq \min \set{v_\pi(a_{ik}b_{kj})}{1 \leq i, j, k \leq m} \\
&\geq& 
\min \set{v_\pi(b_{kj})}{1 \leq  j, k \leq m} = -\ell(b),
\end{eqnarray*}
hence $\ell(ab) \leq \ell(b)$, and similarly we have $\ell(ba) \leq \ell(b)$. If $a \in GL_m(\calO_{k'})$, we have the opposite inequalities and then $\ell(ab) = \ell(b) = \ell(ba)$. 

The statement (2) follows from the fact $det(g) \in \calO_{k'}^\times$ and $K = M_m(\calO_{k'}) \cap G$.
\qed

\vspace{5mm}
\noindent
{\bf 1.2.} 
In this subsection we consider the case $m = 3$ and prove the following proposition. 
\begin{prop}   \label{prop:n=1}
The set $\calR_1^+ \sqcup \calR_1^-$  is a set of complete representatives of $K_1\backslash X_1$, where 
\begin{eqnarray*}
&&
\calR_1^+ = \set{x_\ell = \begin{pmatrix}\pi^\ell & & \\ & 1 & \\ && \pi^{-\ell}\end{pmatrix}}{\ell \geq 0}, \\
&&
\calR_1^- = 
\set{x_{-r} = \begin{pmatrix}\pi^{-r}(1-\eps) &&-\sqrt{\eps}\\ &1 & \\
\sqrt{\eps} & & \pi^r \end{pmatrix}}{1 \leq r \leq e}. 
\end{eqnarray*}
The set $\calR_1^-$ is non-empty only if $e >0$. In that case, for $x \in X_1$, $K_1\cdot x$ has a representative in $\calR_1^-$ if and only if $x \equiv j_3 \pmod{(\pi)}$, and then $r = \ell(x-j_3)$.     
\end{prop}
We write down the group $K_1= K_1^{(od)}$ explicitly for convenience. 
\begin{lem} \label{lem: explicit K1}
\begin{eqnarray*}
&&
K_1 = K_{1,1} \sqcup K_{1,2},\\
&&
K_{1,1}:= \set{g \in B_1j_3B_1 \cap K_1}{g_{31} \in \calO_{k'}^\times} = \set{g \in K_1}{g_{31} \in \calO_{k'}^\times}\\[2mm]
&& \quad = \set{
\begin{pmatrix}\alp &&\\&u&\\&& \alp^{* -1} \end{pmatrix}
\begin{pmatrix}1  &-d^*&f\\&1&d\\&&1 \end{pmatrix}
\begin{pmatrix}  &&1\\&1&-b^*\\ 1&b & c \end{pmatrix}}
{\begin{array}{l}
  \alp \in \calO_{k'}^\times, \; u \in \calO_{k'}^1\\
b, c, d, f \in \calO_{k'}\\
N(b)+c+c^* = N(d)+f+f^* = 0
\end{array}},\nonumber \\[2mm]
&&
K_{1,2}:= \set{g \in K_1}{g_{31} \in (\pi)}\\[2mm]
&&\quad
=\set{
\begin{pmatrix}\alp &&\\& u &\\&&\alp^{* -1} \end{pmatrix}
\begin{pmatrix}1  &&\\ b & 1 &  \\ c& -b^* & 1 \end{pmatrix}
\begin{pmatrix}  1& d & f\\ &1&-d^*\\  &  & 1 \end{pmatrix}}
{\begin{array}{l}
\alp \in \calO_{k'}^\times, \; u \in \calO_{k'}^1\\
b, c \in \pi\calO_{k'}, \; d, f \in \calO_{k'}\\
N(b)+c+c^* = N(d)+f+f^* = 0
\end{array} }.
\end{eqnarray*}
\end{lem}


\bigskip
{\it Proof of }Proposition~\ref{prop:n=1}.  The strategy is similar to [HK-II]-\S 1.2.
We take an element $x \in X_1$, write it as
\begin{eqnarray}   \label{any x}
x = \begin{pmatrix}
a & b & c\\
b^* & d & f\\
c^* & f^* & g
\end{pmatrix}, \qquad a, d, g \in k, \quad b, c, f \in k',
\end{eqnarray}
and show that the orbit $K_1 \cdot x$ has an element $x_\ell$ with $\ell \geq -e$ as in the statement.

\medskip
\noindent
By the fact $x \in G$ and $\Phi_{xj_3}(t) = (t^2-1)(t-1)$, we obtain the following equations
\begin{subequations}
\label{eq:fund_eqs}
\begin{align}
\label{eq:fund_eqs1}
    ag+bf+c^2&=1, \\
\label{eq:fund_eqs2}
    af^*+b(c+d)&=0,\\
\label{eq:fund_eqs3}
    a(c+c^*)+bb^*&=0,\\
\label{eq:fund_eqs4}
    b^*g+(c+d)f&=0,\\
\label{eq:fund_eqs5}
    bf+b^*f^*+d^2&=1,\\
\label{eq:fund_eqs6}
    (c+c^*)g+ff^*&=0,
\end{align}
\end{subequations}
and
\begin{eqnarray} \label{eq:char_poly}
\lefteqn{(t^2-1)(t-1)}\\
&=&(t-c)(t-c^*)(t-d)-(t-c)b^*f^*
-(t-c^*)bf-(t-d)ag-aff^*-bb^*g. \nonumber
\end{eqnarray}

\noindent
We recall that $N(\calO_{k'}^\times) = \calO_k^\times$ and $Tr(\calO_{k'}) = \calO_k$.

\medskip
\noindent
[Case 1] $a \ne 0$ and $v_\pi(a) \leq v_\pi(b)$, or $g \ne 0$ and $v_\pi(g) \leq v_\pi(f)$.

\noindent
By the action of $j_3$, it suffices to consider the case $a \ne 0$. By the action of an element of $K_1$ of type
\begin{eqnarray} \label{eq:basic}
\begin{pmatrix}
1 & 0& 0\\
\lam & 1 &0\\
\mu & -\lam^* & 1
\end{pmatrix} \in K_1, \qquad \lam, \nu \in k' \; \mbox{such that} \; a\lam + b^* = 0, \; N(\lam) + \mu + \mu^* = 0,
\end{eqnarray}
we may assume $b = 0$ in (\ref{any x}). Then, by \eqref{eq:fund_eqs} and \eqref{eq:char_poly}, we have
$$
f = 0, \quad ag + c^2 = 1, \quad c + c^* = 0, \quad d = 1.
$$
Thus, we may assume  
\begin{equation}    \label{eq:Case1}
  x=
  \begin{pmatrix}
    a & 0 & -c_1\sqrt{\epsilon} \\
    0 & 1 & 0 \\
    c_1\sqrt{\epsilon} & 0 & g
  \end{pmatrix}, \quad \begin{array}{l}
a,g,c_1\in k, \; ag+c_1^2\epsilon=1,\\
ag \ne 0, \; v_\pi(a) \geq v_\pi(g), 
\end{array}
\end{equation}
where $g \ne 0$ follows from $\eps \notin k^{\times 2}$. We consider two cases according to $v_\pi(g)$. 
 
\medskip
\noindent
(i) If $v_\pi(g) \leq 0$, we may assume $g = \pi^{-\ell}, \; \ell \geq 0$. 
Then $(\pi^\ell c_1)^2 = \pi^{2\ell} - ag^{-1} \in \calO_k$, and  
\begin{eqnarray*}
K_1 \cdot x &\ni& \begin{pmatrix}
1 & 0 & -\pi^\ell c_1\sqrt{\eps}\\
0& 1 & 0\\
0 & 0 & 1
\end{pmatrix} \cdot 
\begin{pmatrix}
a & 0 & c_1\sqrt{\eps}\\
0& 1 & 0\\
-c_1\sqrt{\eps} & 0 & \pi^{-\ell}
\end{pmatrix} \\
&=&
\begin{pmatrix}
\pi^\ell & 0 & 0 \\
0 & 1 & 0\\
0 & 0 & \pi^{-\ell}
\end{pmatrix}.
\end{eqnarray*}

\noindent
(ii) Assume that $v_\pi(g) > 0$. Since an element of $k$ square modulo $4\pi\calO_k$ is a square in $k$,  the relation $c_1^2\eps \equiv 1 (\pi^{v_\pi(ag)}\calO_k)$ yields $e > 0$ and $v_\pi(g) \leq e$. Further $c_1 \equiv 1 \pmod{(\pi)}$, since $\eps \in 1+4\calO_k$. Thus we see
\begin{eqnarray}
K_1 \cdot x \ni x' = 
\begin{pmatrix}
a & 0 & -(1+b)\sqrt{\eps}\\
0 & 1 & 0\\
(1+b)\sqrt{\eps} & 0 & \pi^r
\end{pmatrix},\quad 
\begin{array}{l}
1 \leq r \leq e,\\
a, b \in \calO_k, v_\pi(a) \geq r,\\
v_\pi(b) < r \; \mbox{or}\; b = 0,
%
\end{array}
\end{eqnarray}
By the equation
\begin{eqnarray*}
(t-1)(t^2-1) = \Phi_{x'j_3}(t) &=&(t-1) \left\{t^2-(1+b)^2\eps - \pi^r a \right\}
\end{eqnarray*}
and $\eps \in 1 +4\calO_k = 1 + (\pi^{2e})$, we have
\begin{eqnarray*}
2b + b^2 + \pi^r a \in \pi^{2e}\calO_k,
\end{eqnarray*}
which yields $b=0$ and $a = \pi^{-r}(1-\eps)$. Hence
\begin{eqnarray}
x' = \begin{pmatrix}
\pi^{-r}(1-\eps) & 0 & -\sqrt{\eps}\\
0 & 1 & 0\\
\sqrt{\eps} & 0 & \pi^r
\end{pmatrix} \; ( = x_{-r}, \; \mbox{say}),\quad  1 \leq r \leq e. 
\end{eqnarray}   
We note $x_{-r}$ appears only if $e > 0$ and cannot be diagonalized by Proposition~\ref{prop:j-type}.

\medskip
\noindent
[Case 2] $ag = 0$. We may assume $a = 0$.

\noindent
By \eqref{eq:fund_eqs} and \eqref{eq:char_poly}, we see $b = 0, \; c = 1, \; d = -1$, and $2g+ff^*=0$. Hence we see 
\begin{eqnarray}   \label{eq:Case2-reduced}
K_1 \cdot x \ni x' = 
\begin{pmatrix}
0 & 0 & 1\\
0 & -1 & h\\
1 & h & -\frac{1}{2}h^2
\end{pmatrix},\quad  h = 0 \mbox{ or } \pi^\ell, \; (\ell \in \Z). 
\end{eqnarray}

\noindent
If $v_\pi(h) \leq e$, then $x'$ satisfies the assumption of Case1, and we have done.
If $h \ne 0$ and $v_\pi(h) > e$, then $\frac{1}{2}h \in \calO_k$. Taking $c \in \calO_{k'}$ such that $\frac{h^2}{4}+c+c^* = 0$, one has 
\begin{eqnarray*}
K_1 \cdot x &\ni& 
\begin{pmatrix}
1 &0&0\\
-\frac{h}{2} &1&0\\
c & \frac{h}{2} & 1
\end{pmatrix} \cdot x' = 
\begin{pmatrix}
0 &0&1\\
0 & -1 & 0\\
1 & 0 & 0
\end{pmatrix}.
\end{eqnarray*}
Hence we have only to consider $h = 0$ in \eqref{eq:Case2-reduced}. i.e. 
\begin{eqnarray} \label{eq:Case2main}
x' = \begin{pmatrix} & & 1\\& -1 & \\1 && \end{pmatrix} = j_3 - 2 y y^*, \quad y = \begin{pmatrix}0\\1\\0\end{pmatrix}.
\end{eqnarray}
For an element $k$ in $K_1$ given by
\begin{eqnarray*}
k &=& \begin{pmatrix}1 & -1 & \frac{-1+\sqrt{\eps}}{2} \\ & 1 & 1\\ &&1 \end{pmatrix}
\begin{pmatrix}&&1 \\ &1 &1 \\ 1 & -1 & \frac{-1+\sqrt{\eps}}{2}\end{pmatrix}\\
&=&
\begin{pmatrix}\frac{-1+\sqrt{\eps}}{2} & \frac{-1-\sqrt{\eps}}{2} & \frac{(-1+\sqrt{\eps})^2}{4} \\ 
1&0 &\frac{1+\sqrt{\eps}}{2} \\ 
1 & -1 & \frac{-1+\sqrt{\eps}}{2}\end{pmatrix}, 
\end{eqnarray*}
one has
\begin{eqnarray*}
k\cdot x' &=& j_3 - k\cdot (2yy^*) = j_3 - 2 (ky)(ky)^*\\
&=&
j_3 - 2 \begin{pmatrix} \frac{-1-\sqrt{\eps}}{2}\\0\\-1\end{pmatrix} 
\begin{pmatrix}\frac{-1+\sqrt{\eps}}{2} &0 &-1\end{pmatrix} \\
&=&
\begin{pmatrix}
\frac{\eps-1}{2} & 0 & -\sqrt{\eps}\\
0 & 1 & 0\\
\sqrt{\eps} & 0 & -2
\end{pmatrix},
\end{eqnarray*}
which is $K_1$-equivalent to 
\begin{eqnarray*}
x_{-e} = \begin{pmatrix}
\pi^{-e}(1-\eps) & 0 & -\sqrt{\eps}\\
0 & 1 & 0\\
\sqrt{\eps} & 0 & \pi^e
\end{pmatrix}.
\end{eqnarray*}
If $e>0$, the above $x_{-e}$ cannot be diagonalized by Proposition~\ref{prop:j-type}  and one of the required representatives. If $e = 0$,  $x_{-0}$ reduces to Case 1-(i), and actually one has
\begin{eqnarray*}
K_1 \cdot x_{-0} \ni \begin{pmatrix}
1 & 0 & \sqrt{\eps}\\
0 & 1 & 0\\
0 & 0 & 1 \end{pmatrix} \cdot 
\begin{pmatrix}
(1-\eps) & 0 & -\sqrt{\eps}\\
0 & 1 & 0\\
\sqrt{\eps} & 0 & 1
\end{pmatrix} = 1_3.
\end{eqnarray*}

\medskip
\noindent
[Case 3] $ag \ne 0, \;  v_\pi(a) > v_\pi(b)$ and $v_\pi(g) > v_\pi(f)$.
We may assume 
\begin{eqnarray} \label{eq:a vs g}
v_\pi(a) \geq v_\pi(g).
\end{eqnarray}
By \eqref{eq:fund_eqs2} and
\eqref{eq:fund_eqs4}, we have
$aN(f)=N(b)g$, which implies
\begin{equation}
  \label{eq:cond_agbf}
  v_\pi(a)-v_\pi(g)=2(v_\pi(b)-v_\pi(f))\geq 0.
\end{equation}
We show $v_\pi(c) \geq 0$. If it was not, by \eqref{eq:fund_eqs1} and the assumption of Case 3, we had
\begin{eqnarray*}
v_\pi(bf) = 2v_\pi(c) \leq -2, \quad \mbox{hence \; } v_\pi(f) < 0.
\end{eqnarray*}
Then
\begin{eqnarray*}
(cf^{-1})^2 = f^{-2}-(af^{-1})(gf^{-1}) - bf^{-1} \in \calO_{k'} \quad (\mbox{by } \eqref{eq:fund_eqs1}),
\end{eqnarray*}
hence $cf^{-1}, c^*f^{-1} \in \calO_{k'}$ and so
\begin{eqnarray*}
(\pi) \ni (cf^{-1}+c^* f^{-1})gf^{* -1} = (c+c^*)gN(f)^{-1} = -1 \quad (\mbox{by } \eqref{eq:fund_eqs6}), 
\end{eqnarray*}
which is a contradiction. Hence $v_\pi(c) \geq 0$, and again by \eqref{eq:fund_eqs6} and the assumption of Case 3, we have
\begin{eqnarray*}
v_\pi(f) < v_\pi(g) \leq 2v_\pi(f), \quad \mbox{hence } v_\pi(f) > 0.
\end{eqnarray*}
Then, by \eqref{eq:fund_eqs5}, $d \in \calO_k^\times$ and $x \in M_{2n+1}(\calO_{k'}) \cap X$, hence $\ell(x) = 0$. 
Set $r = v_\pi(f) > 0$.
Then 
\begin{eqnarray*}
&&
c+c^*+d = 1, \quad (\mbox{by } \eqref{eq:char_poly})\\ 
&&
c +d = c^* + d \equiv 0 \pmod{(\pi^{r+1})}, \quad (\mbox{by } \eqref{eq:fund_eqs4}) 
\end{eqnarray*}
hence
\begin{eqnarray*}
c \equiv c^* \equiv 1 \pmod{(\pi^{r+1})}, \quad c+c^* \equiv 2 \pmod{(\pi^{r+1})}.
\end{eqnarray*}
Then by \eqref{eq:fund_eqs6}
\begin{eqnarray*}
2r = v_\pi(c+c^*)+v_\pi(g) \geq \min\{e, r+1\}+r+1,
\end{eqnarray*}
hence we see
\begin{eqnarray*}
r > e, \quad v_\pi(c+c^*) = e, \quad v_\pi(g) = 2r-e > r.
\end{eqnarray*}
Setting $v_\pi(b) = m (\geq r)$, one has $v_\pi(a) = 2m-e$ by \eqref{eq:fund_eqs3}.
We may take the unit part of $a$ off, and assume $x$ becomes
\begin{eqnarray} \label{eq:Case3}
x = \begin{pmatrix}
\pi^{2m-e} & \pi^m u & c \\
\pi^m u^* & d & \pi^rv\\
c^* & \pi^rv^* & \pi^{2r-e}w
\end{pmatrix}, \quad
\begin{array}{l}
m \geq r > e\\
c, u, v \in \calO_{k'}^\times, \; d, w \in \calO_k^\times\\
c \equiv c^* \equiv 1 \mod(\pi^{r+1}), \; d \equiv -1 \mod(\pi^{r+1}),\\
c + c^* + d = 1,
\end{array}
\end{eqnarray}
and the set of equations \eqref{eq:fund_eqs} becomes
\begin{eqnarray*}
  &&
\pi^{2(m+r-e)}w+\pi^{m+r}uv+c^2=1,\\
  &&
c+d=-\pi^{m+r-e}u^{-1}v^*,\\
  &&
c+c^*=-\pi^euu^*,\\
  &&
c+d=-\pi^{m+r-e}u^*v^{-1}w,\\
  &&
\pi^{m+r}(uv+u^*v^*)+d^2=1,\\
  &&
c+c^*=-\pi^e vv^*w^{-1}.
\end{eqnarray*}
Then together with \eqref{eq:Case3}, we have
\begin{eqnarray*}
&&
w = (uu^*)^{-1}vv^*,\quad c+c^* = -\pi^euu^*,\\
&&
d = 1 - (c+c^*) = 1 + \pi^euu^*,\\
&&
c = 1-(c^*+d) = 1+ \pi^{m+r-e}u^{* -1}v.\\
\end{eqnarray*}
Now we have
\begin{eqnarray*}
x' := \begin{pmatrix}
 u^{-1} & & \\ & 1 & \\ &&u^*
\end{pmatrix} \cdot x 
&=& 
\begin{pmatrix}
\pi^{2m-e}(uu^*)^{-1} & \pi^m & c\\
\pi^m & d & \pi^ruv\\
c^* & \pi^ru^*v^* & \pi^{2r-e}vv^*
\end{pmatrix} \\[2mm]
&=&
j_3 + \frac{\pi^e}{uu^*} \begin{pmatrix}
\pi^{2(m-e)} & \pi^{m-e}uu^* & \pi^{m+r-2e}uv\\
\pi^{m-e}uu^* & (uu^*)^2 & \pi^{r-e}u^2u^*v\\
\pi^{m+r-2e}u^*v^* & \pi^{r-e}uu^{* 2}v^* & \pi^{2(r-e)}uu^*vv^*
\end{pmatrix} \\[2mm]
&=&
j_3 + \frac{\pi^e}{uu^*} yy^*, \quad y = \begin{pmatrix} \pi^{m-e}\\ uu^* \\ \pi^{r-e}u^*v^*\end{pmatrix}.
\end{eqnarray*}
We will show there exists some $k \in K_1$ for which the second entry $(ky)_2$ of $ky$ =0.  Then $k\cdot x'$ has the shape $\begin{pmatrix} \cdot & 0 & \cdot\\ 0 & 1 & 0\\ \cdot& 0 &\cdot \end{pmatrix}$, which reduces to Case 1 or Case 2 (if $e >0$,  $x$ is equivalent to some $x_{-\ell}, \; 1 \leq \ell \leq e$,  by Proposition~\ref{prop:j-type}). 
Set
\begin{eqnarray} \label{eq:to find k}
k = 
\begin{pmatrix}
1 &-\alp^* & \beta\\
0 & 1 & \alp\\
0 & 0 & 1
\end{pmatrix}
\begin{pmatrix}
0 &0 & 1\\
0 & 1 & -\gamma^* \\
1 & \gamma & \delta
\end{pmatrix} = \begin{pmatrix} &\cdot &\\\alp & 1+\alp\gamma & \alp\delta - \gamma^* \\ &\cdot &\end{pmatrix}, 
\end{eqnarray}
and solve
\begin{eqnarray} \label{eq:target 0}
\big((ky)_2 = \big)\,  \alp \pi^{m-e} + (1+\alp\gamma)uu^* + (\alp\delta-\gamma^*)\pi^{r-e}u^*v^* = 0
\end{eqnarray}
under the condition $\alp, \beta, \gamma, \delta \in \calO_{k'}$ and $N(\alp)+\beta + \beta^* = N(\gamma) + \delta + \delta^* = 0$, which is equivalent that $k$ of \eqref{eq:to find k} becomes an element in $K_1$.

If \eqref{eq:target 0} is satisfied, we see $1+\alp\gamma \in (\pi^{r-e}) \subset (\pi)$, and $\alp, \gamma \in \calO_{k'}^\times$. Writing $1 + \alp\gamma = \pi^{r-e}\lam$ with $\lam \in \calO_{k'}$, we have 
$\delta = \alp^{-1}\gamma^* -\pi^{m-r}(u^*v^*)^{-1}-\alp^{-1}\lam uv^{* -1}$ by \eqref{eq:target 0}. Then, since $\alp \ne 0$,  the condition $N(\gamma)+\delta+\delta^*=0$ is equivalent to 
\begin{eqnarray} \label{eq:(1)}
N(\alp)N(\gamma) +  \alp\gamma + \alp^*\gamma^* - \pi^{m-r}N(\alp)((uv)^{-1}+(u^*v^*)^{-1}) - \alp^*\lam uv^{* -1} - \alp\lam u^*v^{-1} = 0, \nonumber
\end{eqnarray}
then, since $1+\alp\gamma = \pi^{r-e}\lam$, it becomes
\begin{eqnarray} \label{eq:(2)}
\pi^{2(r-e)}N(\lam)  - \pi^{m-r}N(\alp)((uv)^{-1}+(u^*v^*)^{-1}) - \alp^*\lam uv^{* -1} - \alp\lam^* u^*v^{-1}= 1. 
\end{eqnarray}  
Setting $\lam = \alp u^{-1}v^*\mu$ with $\mu \in \calO_{k'}$, \eqref{eq:(2)} is equivalent to 
\begin{eqnarray} \label{eq:(3)}
N(\alp)\left(\pi^{2(r-e)}N(u)^{-1}N(v)N(\mu)  - \pi^{m-r}((uv)^{-1}+(u^*v^*)^{-1}) - (\mu + \mu^*)\right) = 1.
\end{eqnarray}  
Since $r >e$, we may choose $\mu \in \calO_k$ for which the latter factor in the left hand side of \eqref{eq:(3)} becomes a unit in $\calO_k$, then for suitable $\alp \in \calO_{k'}^\times$ we establish the identity \eqref{eq:(3)}. Finally taking $\beta$ such as $N(\alp)+\beta+\beta^*=0$, we obtain $k \in K_1$ for which $(ky)_2=0$, which establishes \eqref{eq:target 0}. 

\medskip
\noindent
Thus we have shown that any $K_1$-orbit in $X_1$ has a representative in $\calR_1^+ \sqcup \calR_1^-$ and $K_1 \cdot x \cap \calR_1^- \ne \emptyset$ if and only if $e > 0$ and $x \equiv j_3 \pmod{(\pi)}$. It is known that each $x_\ell, \; \ell \geq 0$ gives a different $GL_3(\calO_{k'})$-orbit in $\set{x \in GL_3(k')}{x^* = x}$,  hence it gives a different $K_1$-orbit in $X_1$. 
For $x_{-r}, \; 1 \leq r \leq e$,  $r = \ell(x_{-r}-j_3) $ is an invariant of $K_1\cdot x_{-r}$, since $h\cdot (x_{-r}-j_3) = h \cdot x -j_3$ for any $h \in K_1$. Hence $\calR_1^+ \sqcup \calR_1^-$ forms a set of complete representatives of $K_1 \backslash X_1$.
\qed  

\vspace{5mm}
\noindent
{\bf 1.3.} 
In this subsection we will show Theorem~\ref{th:Cartan}-(1) for the case $m = 2n+1$ with $n \geq 2$. Our strategy is similar to HK-II (we have to be careful for $e > 0$.)


\begin{lem}  \label{lem:kakunin}
Let $n \geq 2$. Then every $x \in X_n$ has a minimal entry except $(n+1, n+1)$-entry.
\end{lem}

\proof
Assume $x$ had unique minimal entry $\pi^{-\ell}u, \; u \in \calO_k^\times$ at $(n+1, n+1)$ and denote by $E'$ the $(n+1,n+1)$-matrix unit. Then $\pi^\ell x \in M_{2n+1}(\calO_{k'})$ and  
$$
\pi^{2\ell}1_{2n+1} = \pi^\ell x \cdot j_{2n+1} \pi^\ell x j_{2n+1} \equiv  u^2E' \mod (\pi),
$$
which is impossible.
\qed

\bigskip
The following lemmas can be shown in the same way as in [HK-II], so omit the details of proofs. 

\begin{lem}  \label{lem:A1-n}
Let $n \geq 2$ and assume that 

$(A1,n): $  $x \in X_n$ has a minimal entry in the diagonal except the $(n+1, n+1)$-entry. 

\noindent
Then $K\cdot x$ contains a hermitian matrix of the type
\begin{eqnarray*}
\left( \begin{array}{c|c|c}
\pi^\ell & 0 & 0\\ 
\hline
0 & y & 0\\  
\hline
0 & 0 & \pi^{-\ell}
\end{array} \right), \qquad y \in X_{n-1}, \; \ell = \ell(x) \geq \ell(y).
\end{eqnarray*}
\end{lem}

{\it Outline of a proof.}
By the action of $W$, we may assume that $(2n+1, 2n+1)$-entry of $x$ is minimal, and then we may arrange it into the above shape by the suitable $K$-action.
\qed

\begin{lem}   \label{lem:A2, A3}
Let $n \geq 2$ and assume one of the following conditions hold: 

$(A2,n):$ $x \in X_n$ has a minimal entry outside of the diagonal, the anti-diagonal, the $(n+1)$th row, and the $(n+1)$th column. 

$(A3,n):$  $x$ has a minimal entry and a non-minimal entry in the anti-diagonal except the $(n+1,n+1)$-entry.  

\medskip
\noindent
Then $K \cdot x$ contains a hermitian matrix of type
$$
\left( \begin{array}{c|c|c}
\begin{array}{cc}
\pi^{\ell}&0\\ 0 &\pi^{\ell}  
\end{array}
 & 0 & 0\\
\hline 
0 & y & 0\\
\hline  
0 & 0 & 
\begin{array}{cc}
\pi^{-\ell}&0\\ 0 &\pi^{-\ell}  
\end{array}
\end{array}\right), \quad 
\begin{array}{l}
y \in X_{n-2}, \, \ell=\ell(x) \geq \ell(y) \mbox{\; if\; } n \geq 3,\\
y = 1 \mbox{\; if\; } n = 2.
\end{array} 
$$
\end{lem}

{\it Outline of a proof.}
Write $\ell = \ell(x)$. Under the condition $(A2, n)$, such minimal entries appear in pair, since $x$ is hermitian. Then, through the action of $W$ and $GL_2(\calO_{k'})$, we may assume that the lower right $2$ by $2$ block of $x$ is $\twomatrix{\pi^{-\ell}}{0}{0}{\pi^{-\ell}}$, then we may arrange it into the above shape by the suitable $K$-action. As for $y$, it is clear that $y = y^*$ and $\ell(y) \leq \ell$. Considering the characteristic function of $xj_{2n+1}$, we see $y \in X_{n-2}$ if $n \geq 3$ and $y = 1$ if $n = 2$.

As for the condition $(A3,n)$, together with the action of $W$, we may assume 
the upper right $2$ by $2$ block of $x$ is $\twomatrix{a}{\xi}{b}{c}$ such that $v_\pi(\xi) = -\ell$ and $a, b, c \in \pi^{-\ell+1}\calO_{k'}$.  Then, by the action of
$$
h = 
\left( \begin{array}{c|c|c}
\begin{array}{cc}
1&1\\ 0 &1  
\end{array}
 & 0 & 0\\
\hline 
0 & 1_{2n-3} & 0\\
\hline  
0 & 0 & 
\begin{array}{cc}
1&-1\\ 0 &1  
\end{array} 
\end{array}\right),
$$
the $(1,2n)$-entry becomes $a+b-c-\xi$, and its $v_\pi$-value is $\ell = \ell(x)$. Thus it reduces to the case $(A2,n)$. 
\qed  

\bigskip
We have to consider the remaining case, which is the assumption $(A4, n)$ below. 
The statement for the non-dyadic case is a refinement of [HK-II, Lemma~1.9].

\begin{lem}
\label{lem:hardest}
  Let $x \in X_n$ with $n \geq 2$. Assume that 

$(A4,n):$  any minimal entry of $x \in X_n$ stands in the anti-diagonal,  the $(n+1)$th row, or the $(n+1)$th column, 
and that all the anti-diagonal entries except $(n+1, n+1)$ are minimal.
 
\noindent
Then $\ell(x) = 0$, and  \\
{\rm (i)} if $k$ is non-dyadic, then $K \cdot x$ contains $1_{2n+1}$;\\ 
{\rm (ii)} if $k$ is dyadic, then $K\cdot x$ contains 
\begin{eqnarray} \label{eq:except n}
\begin{pmatrix}
    \pi^{\mu_1}(1-\eps) &&&&&& -\sqrt{\eps} \\
    & \ddots &&&& \iddots & \\
    &&\pi^{\mu_n}(1-\eps)&& -\sqrt{\eps}&& \\
&& & 1 &&&\\    
&&\sqrt{\eps}&&\pi^{-\mu_n} && \\
    & \iddots &&&& \ddots & \\
    \sqrt{\eps} &&&&&& \pi^{-\mu_1}
  \end{pmatrix}, 
\end{eqnarray}
where $0>\mu_1 \geq \mu_2 \geq \cdots \geq \mu_n \geq -e$.

\end{lem}

\proof
  By the action of $W$ and suitable $K$-action, we may assume $x$ has the shape 
  \begin{equation} \label{eq:stage1}
    x=\left(
    \begin{array}{cccc|c|cccc}
      &&&& \alpha &&&& \xi_1 \\
      &*&&&    0   &&*&\iddots& \\
      &&&& \vdots &&\iddots&*& \\
      &&&& 0 & \xi_n &&& \\
      \hline
      \alpha^* & 0 &\cdots  &0 & u&0 &\cdots& 0 &\kappa^*\\
      \hline
      &&& \xi_n^*& 0 &&&& \\
      &*&\iddots&& \vdots &&&& \\
      &\iddots&*&& 0 &&&*& \\
       \xi_1^*&&&& \kappa &&&&\\
    \end{array}
  \right),
  \end{equation}
  where $\alp = 0$ or $v_\pi(\alpha)\leq v_\pi(\kappa)$. Write $\ell = \ell(x)$.
Since $xj_{2n+1}xj_{2n+1}=1_{2n+1}$, by its $(1,2n+1)$-entry and $(2n+1,1)$-entry, we have
\begin{eqnarray*}
\alpha\alpha^*\equiv \kappa\kappa^*\equiv0 \pmod{(\pi^{-2\ell+1})},
\end{eqnarray*}
hence $\alp$ nor $\kappa$ is not minimal, which means any entry outside of the anti-diagonal is non-minimal.  
Setting $u_i = \pi^\ell \xi_i \in \calO_{k'}^\times, \; 1 \leq i \leq n$ and $u_0 = \pi^\ell u, \; u \in \calO_k$, we have 
\begin{eqnarray*}
\pi^{2\ell}1_{2n+1} &=& \pi^{2\ell}xj_{2n+1}xj_{2n+1} = (\pi^\ell x) j_{2n+1}(\pi^\ell x) j_{2n+1} \\
&\equiv&  
Diag(u_1^2, \ldots, u_n^2, u_0^2, u_n^{* 2}, \ldots, u_1^{* 2}) \; \pmod{(\pi)},
\end{eqnarray*}
hence 
\begin{eqnarray*} \label{eq:xi-u}
\ell = 0, \quad \xi_i^2 \equiv u^2 \equiv 1 \pmod{(\pi)}, 
\end{eqnarray*}
and
\begin{eqnarray*} 
&& \label{eq:stage1a}   
x\equiv \left(
    \begin{array}{ccc|c|ccc}
      &&& &&& \xi_1 \\
      &&& &  &\iddots& \\
      &&&  & \xi_n && \\
      \hline
       & & & u& && \\
      \hline
      && \xi_n&&&& \\
      &\iddots&&&&& \\
       \xi_1&&&&&&\\
    \end{array}
  \right) \pmod{(\pi)}, \\[2mm]
&&
\xi_i \equiv \xi_i^* \equiv \pm 1\pmod{(\pi)}, \quad u \equiv \pm 1\pmod{(\pi)}.\nonumber
\end{eqnarray*}

First we consider the non-dyadic case, where $1 \not\equiv -1 \pmod{(\pi)}$. By the characteristic polynomial, we have 
\begin{eqnarray*}
&&
(t^2-1)^n(t+1) \equiv (t-u)\prod_{i=1}^n (t-\xi_i)^2 \pmod{(\pi)}, \quad n \geq 2,
\end{eqnarray*}
hence we may assume $\xi_1 \not\equiv \xi_2 \pmod{(\pi)}$, after suitable $W$-action if necessary.
For
\begin{eqnarray*}
h = 
\left(\begin{array}{cc|c|cc}
1 &&&&\\
&1&&&\\
\hline
&&1_{2n-3}&&\\
\hline
1 &&&1&\\
&-1&&&1 
\end{array} \right) \in K,
\end{eqnarray*}
we have
\begin{eqnarray*}
h \cdot x \equiv \left(\begin{array}{cc|c|cc}
&&&&\xi_1\\
&&&\xi_2&\\
\hline
&&C&&\\
\hline
&\xi_2&&&\xi_1-\xi_2\\
\xi_1&&&\xi_1-\xi_2& 
\end{array} \right) \pmod{(\pi)}, 
\end{eqnarray*}
where $Cj_{2n-3} = Diag(\xi_3, \ldots, \xi_n,u,\xi_n,\ldots,\xi_3)$ and $\xi_1-\xi_2 \not\equiv 0 \pmod{(\pi)}$.
Since $h\cdot x$ satisfies $(A2, n)$ with $\ell = 0$, we see 
\begin{eqnarray*}
K \cdot x \ni 
\left( \begin{array}{c|c|c}
1_2
 & 0 & 0\\
\hline 
0 & y & 0\\
\hline  
0 & 0 & 1_2
\end{array}\right), \quad y \in X_{n-2}, \; \ell(y) = 0, \; \mbox{or\; } y = 1.
\end{eqnarray*} 
Then,  we see inductively $K\cdot x$ contains $1_{2n+1}$.

\bigskip
Now we assume $k$ is dyadic. Then $x \equiv j_{2n+1} \; (\mod (\pi))$ and $x$ is not diagonalizable, by Proposition~\ref{prop:j-type}.   

\medskip
\noindent
[Case 1] We assume that there is a non-zero entry of $x$ outside of  the anti-diagonal, $(n+1)$-th row and $(n+1)$-th column, i.e. within $*$-places of \eqref{eq:stage1}. Let $\ell$ be the minimal $v_\pi$-value within those entries. Then, after suitable $K$-action, we may assume the $(2n+1, 2n+1)$-entry of $x$ is $\pi^\ell$.   
Further, after suitable $K$-action, we may assume 
  \begin{equation*} \label{eq:stage2}
    x=\left(
    \begin{array}{cccc|c|cccc}
      c&&&& \alpha &&&& \xi_1 \\
      &*&&&    0   &&*&\iddots&0 \\
      &&&& \vdots &&\iddots&*&\vdots \\
      &&&& 0 & \xi_n &&&0 \\
      \hline
      \alpha^* & 0 &\cdots  &0 & u&0 &\cdots& 0 &\kappa^*\\
      \hline
      &&& \xi_n^*& 0 &&&&0 \\
      &*&\iddots&& \vdots &&&&\vdots \\
      &\iddots&*&& 0 &&&*&0 \\
       \xi_1^*&0&\cdots&0& \kappa &0&\cdots&0&\pi^\ell\\
    \end{array}
  \right),
  \end{equation*}
  where $v_\pi(c) \leq \ell$, ``$\alp = 0$ or $v_\pi(\alpha)\leq v_\pi(\kappa)$", and ``$\kappa = 0$ or $v_\pi(\kappa) < \ell$". 
By $(i, 1)$-entry with $i \ne n+1$ of $xj_{2n+1}xj_{2n+1} = 1_{2n+1}$ and hermitianess of $x$, we have 
  \begin{equation*} \label{eq:stage3}
    x=\left(
    \begin{array}{cccc|c|cccc}
      c&0&\cdots&0& \alpha &0&\cdots&0& \xi_1 \\
      0&&&&    0   &&*&\iddots&0 \\
      \vdots&&*&& \vdots &&\iddots&*&\vdots \\
      0&&&& 0 & \xi_n &&&0 \\
      \hline
      \alpha^* & 0 &\cdots  &0 & u&0 &\cdots& 0 &\kappa^*\\
      \hline
      0&&& \xi_n^*& 0 &&&&0 \\
      \vdots&*&\iddots&& \vdots &&*&&\vdots \\
      0&\iddots&*&& 0 &&&&0 \\
       \xi_1^*&0&\cdots&0& \kappa &0&\cdots&0&\pi^\ell\\
    \end{array}
  \right).
  \end{equation*}
From the above $x$, we set
\begin{eqnarray}   \label{eq:9points z}
z  = \begin{pmatrix} c &\alp & \xi_1\\ \alp^* & u & \kappa^*\\ \xi_1^* & \kappa & \pi^\ell,
\end{pmatrix} 
\end{eqnarray}
and see $z$ or $-z$ is an element of $X_1$ with $\ell(z) = 0$ (cf. \eqref{eq:m=2,3}).
By the action of $K_1 = U(j_3)(\calO_{k'})$ through the embedding
\begin{eqnarray*} \label{eq:embedding1}
\begin{array}{cll}
K_1  &\longrightarrow &K = K_n,\\
h = (h_{ij}) & \longmapsto &
\wt{h} = \begin{pmatrix}
h_{11}&&h_{12} && h_{13}\\
& 1_{n-1}& &&\\
h_{21} && h_{22} && h_{23}\\
&&& 1_{n-1} &\\
h_{31} && h_{32} && h_{33}
\end{pmatrix}
\end{array},
\end{eqnarray*} 
we may change $z$ in $x$ of \eqref{eq:9points z} into $\pm x_{-r}^{(1)}$ for some $r$ with $1 \leq r \leq e$, where the superscript $(1)$ indicates the size $(m,n) = (3,1)$.  
When $-x_{-r}^{(1)}$ appears as $z$, by the action of $K_0 := U(j_2)(\calO_{k'})$ (tentative naming) through the embedding
\begin{eqnarray*} \label{eq:embedding0}
\begin{array}{cll}
K_0  &\longrightarrow &K = K_n,\\
h = (h_{ij}) & \longmapsto &
\wt{h} = \begin{pmatrix}
h_{11}&&h_{12}\\
& 1_{2n-1}& \\
h_{21} && h_{22}
\end{pmatrix}
\end{array},
\end{eqnarray*} 
we may  change $z = -x_{-r}^{(-1)}$ of $x$ into 
\begin{eqnarray*}
\begin{pmatrix}
 \pi^{-r_1}(1-\eps) &&-\sqrt{\eps}\\
&-1&\\
\sqrt{\eps} && \pi^{r_1}
\end{pmatrix}, 
\end{eqnarray*}
where $r_1$ might be changed from $r$ but still $1 \leq r_1 \leq e$.
Anyway, we see 
\begin{eqnarray*}
K \cdot x \ni \begin{pmatrix}
\pi^{-r}(1-\eps) & &-\sqrt{\eps}\\
&y&\\
\sqrt{\eps} && \pi^r
\end{pmatrix}, \quad 
\begin{array}{l}
1 \leq r \leq e, \\
y = y^* \in M_{2n-1}(\calO_{k'}), \; y \equiv j_{2n-1} \pmod{(\pi)}.
\end{array}
\end{eqnarray*}
Since $\Phi_{xj_{2n+1}}(t) = (t^2-1)\Phi_{yj_{2n-1}}(t)$, we see $y \in X_{n-1}$, and $y$ satisfies $(A4, n-1)$.  
By an inductive procedure, we see the $K$-orbit of $x$ contains a matrix of type \eqref{eq:except n}. 

\medskip
\noindent
[Case 2] We consider the remaining situation of Case 1, i.e. any entry of $*$-places of $x$ in \eqref{eq:stage1} is $0$. Then $\alp = \kappa = 0$ follows from $xj_{2n+1}xj_{2n+1} = 1_{2n+1}$, and we have 
$$
z_1 = 
\begin{pmatrix}
&&\xi_1\\&u&\\\xi_1&&
\end{pmatrix}, \quad \xi_1 = \xi_1^* = \pm 1, \; u = \pm 1,  
$$
in stead of $z$ of \eqref{eq:9points z}. By the same procedure as Case 1, we see $K\cdot x$ contain a matrix of type \eqref{eq:except n}.
\qed  

\bigskip
By Proposition~\ref{prop:n=1}, Lemma~\ref{lem:A1-n}, Lemma~\ref{lem:A2, A3} and Lemma~\ref{lem:hardest}, we see for every $x \in X$, $K\cdot x$ has a representative of shape $x_\lam$ for some $\lam \in \wt{\Lam_n^+}$, which completes the proof of Theorem~\ref{th:Cartan}-(1).
\qed

\vspace{2cm}
\Section{Spherical function on $X$}
{\bf 2.1.} 
We consider $m = 2n$ or $2n+1$, and write $X = X_n, \; G = G_n, \; B = B_n, \; K = K_n$.   
For $g \in G$, we denote by $d_i(g)$ the determinant of the lower right $i$ by $i$ block of $g$. Then $d_i(x), \; 1 \leq i \leq n$, are relative $B$-invariants on $X$ associated with rational character $\psi_i$ of $B$, where
\begin{eqnarray*} \label{eq:rel inv}
d_i(p\cdot x) = \psi_i(p)d_i(x), \quad \psi_i(p) = N(d_i(p)), \quad (x \in X, \; p \in B).
\end{eqnarray*}
We set
\begin{eqnarray*} \label{eq:X-open}
X^{op} = \set{x \in X}{d_i(x) \ne 0, \; 1 \leq i \leq n}, 
\end{eqnarray*}
then $X^{op}(\ol{k})$ is a Zarisky open $B(\ol{k})$-orbit. 
For $x \in X$ and $s \in \C^n$, we consider the integral 
\begin{eqnarray}  \label{eq:def of sph} 
\omega(x;s) = \int_{K}\abs{\bfd(k\cdot x)}^s dk, \quad 
\abs{\bfd(y)}^s = \left\{\begin{array}{ll}
\prod_{i=1}^n \abs{d_i(y)}^{s_i} & \mbox{if } y \in X^{op}\\
0 & \mbox{otherwise}
\end{array} \right.
\end{eqnarray}
where $\abs{\;}$ is the absolute value on $k$ normalized by $\abs{\pi}= q^{-1}$, $dk$ is the normalized Haar measure on $K$. The integral in \eqref{eq:def of sph} is absolutely convergent if $\real(s_i) \geq 0, \; 1 \leq i \leq n$, and continued to a rational function of $q^{s_1}, \ldots. q^{s_n}$, and we use the notation $\omega(x;s)$ in such sense. We call $\omega(x;s)$ a spherical function on $X$, since it becomes an $\hec$-common eigenfunction on $X$ (cf. \cite[\S 1]{JMSJ}, \cite[\S 1]{French}). Indeed, $\hec$ is a commutative $\C$-algebra spanned by all the characteristic functions of double cosets $KgK, \; g \in G$, by definition, and we see
\begin{eqnarray}
(f*\omega(x;s))(x) &{\big(}& = \int_G\, f(g)\omega(g^{-1}\cdot x; s)dg \big) \nonumber\\
&=&\lam_s(f)\omega(x;s), \quad (f \in \hec),
\end{eqnarray}
where $dg$ is the Haar measure on $G$ normalized by $\int_K dg = 1$, and $\lam_s$ is the $\C$-algebra homomorphism (Satake transform) defined by 
\begin{eqnarray} \label{Satake trans s}
\begin{array}{rll}
\lam_s: \hec & \longrightarrow & \C(q^{s_1}, \ldots, q^{s_n}) ,\\
f & \longmapsto & \displaystyle{\int_{B}} f(p)\abs{\psi(p)}^{-s}\delta(p)dp.
\end{array}
\end{eqnarray}
Here $\abs{\psi(p)}^{-s} = \prod_{i=1}^n \abs{\psi_i(p)}^{-s_i}$, $dp$ is the left-invariant measure on $B$ such that $\int_{B\cap K}dp = 1$, and $\delta(p)$ is the modulus character of $dp \; (d(pq) = \delta(q)^{-1}dp)$. 

It is convenient to introduce a new variable $z$ which is related to $s$ by
\begin{eqnarray}
&&
s_i = -z_i+z_{i+1}-1+\frac{\pi\sqrt{-1}}{\log q} \quad 1 \leq i \leq n-1 \nonumber\\
&& \label{eq:change of var}
s_n= \left\{\begin{array}{ll}
-z_n-\frac12  & \mbox{if } m = 2n\\
-z_n-1 + \frac{\pi\sqrt{-1}}{2\log q} & \mbox{if } m = 2n+1,
\end{array} \right.
\end{eqnarray}
and denote $\omega(x;z) = \omega(x;s)$ and $\lam_z = \lam_s$. 

\begin{rem} \label{rem:shift}
{\rm 
For the case $m = 2n$, we have slightly changed the shift of invariants in \eqref{eq:def of sph} and \eqref{eq:change of var} from those in \cite{HK}, where we set $s_n = -z_n-\frac{1}{2}+\frac{\pi\sqrt{-1}}{\log q}$. By the explicit formula of $\omega(x;z)$  (Theorem~\ref{th:explicit}), we will see the present $\omega(x; z)$ takes the same value as before on $G\cdot x_0$ and the multiple by $(-1)$ on $G\cdot x_1$, and we will explain the reason of this change in Remark~\ref{rem:z0}. 
} 
\end{rem}

Keeping the relation \eqref{eq:change of var}, we see
\begin{eqnarray} \label{eq:psi-s,z} 
&&
\abs{\psi(p)}^s =  \delta^{\frac12}(p) \prod_{i=1}^n\abs{N(p_i)}^{z_i} \times \left(\begin{array}{ll}
1 &\mbox{if } m=2n\\
(-1)^{v_\pi(p_1\cdots p_n)} &\mbox{if }m=2n+1
\end{array} \right), \nonumber\\
&& \label{eq:Satake trans z}
\lam_z: \hec \stackrel{\sim}{\longrightarrow} \C[q^{\pm 2z_1}, \ldots, q^{\pm 2z_n}]^W, 
\end{eqnarray}
where $W$ is the Weyl group of $G$ with respect to the maximal $k$-split torus in $B$. Then $W \cong S_n \ltimes (\pm 1)^n$ acts on $s$ and $z$ through rational characters of $B$, where $W $ is generated by $S_n$ and $\tau$, $S_n$ acts on $z$ by the permutation of indices and $\tau(z) = (z_1,\ldots, z_{n-1},-z_n)$. 
 The functional equation with respect to $S_n$ is reduced to the case of unramified hermitian forms as follows.
Define an embedding $K_0 = GL_n(\calO_{k'}) \longrightarrow K = K_n$ by
\begin{eqnarray*}
K_0 \ni h \longmapsto \wt{h} = \left\{\begin{array}{ll}
\begin{pmatrix}
j_nh^{* -1}j_n  & \\
&h
\end{pmatrix}  & \mbox{if } m = 2n\\
\begin{pmatrix}
j_nh^{* -1}j_n & & \\
&1&\\
&&h
\end{pmatrix} & \mbox{if } m = 2n+1
\end{array} \right\} \in K. 
\end{eqnarray*}
Then, as is considered in \cite{HK} and \cite{HK2}, we have
\begin{eqnarray*}
\omega(x; s) &=& \int_{K_0}dh \int_{K} \abs{\bfd(k\cdot x)}^s dk\\
&=&
\int_{K_0} \int_{K} \abs{\bfd(\wt{h}k\cdot x)}^s dk dh \\
&=&
\int_{K} \zeta_*^{(h)}(D(k\cdot x);s)dk,
\end{eqnarray*}
where $D(k\cdot x)$ is the lower $n$ by $n$ block of $k\cdot x$, and $\zeta_*^{(h)}(y;s)$ is a spherical function on hermitian matrices $\calH_n(k')$. Since the behaviour of $\zeta_*^{(h)}(y; s)$ is independent of the residual characteristic (cf. \cite{JMSJ}), we may quote the the following from \cite[Theorem~2.1]{HK} and \cite[Theorem~2.1]{HK2}.

\begin{prop} The function $G_1(z)\cdot \omega(x; s)$ is invariant under the action of $S_n$ on $z$, where
\begin{eqnarray*}
&&
G_1(z) = \prod_{1 \leq i < j \leq n}\, \frac{1+q^{z_i+z_j}}{1-q^{z_i-z_j-1}}.
\end{eqnarray*}
\end{prop}

\vspace{5mm}
\noindent
{\bf 2.2.} In this subsection we study $\omega(x; s)$ for the case $(m,n)=(2,1)$ and show the following.  
As a result, we see again the set $\xx_1^{ev} = \set{x_\lam}{\lam\in\wt{\Lam_1^+}}$ forms a set of complete representatives of $K_1\backslash X_1$, since $\omega(x_\lam;z)$ takes different value for each $\lam$ for generic $z$. 
\begin{prop} \label{prop: even n=1}
For $x_\lam \in \xx_1^{ev}$, one has 
\begin{eqnarray*}\label{eq:sph m=2}
\omega(x_\lam; z) =  
\frac{q^{-\frac{\lam}{2}} q^{ez} }{ 1+q^{-1}} 
\left(\frac{q^{-(\lam+e)z}(1-q^{2z-1})}{1-q^{2z}} +  \frac{q^{(\lam+e)z}(1-q^{-2z-1})}{1-q^{-2z}} \right) .
\end{eqnarray*}
In particular, for any $x \in X_1^{(ev)}$, 
 \begin{eqnarray} \label{eq:fun-eq m=2}
q^{-ez}\omega(x;z) \in \C[q^z+q^{-z}], \quad \omega(x; z) =  q^{2ez}\,\omega(x; -z).
\end{eqnarray}
\end{prop}

\medskip
We have proved for $e=0$ in \cite[Proposition~2.4]{HK}, but we give a unified proof for $e \geq 0$ here. 
%
It is easy to see 
\begin{eqnarray*}
&&
K_1 = K_{1,1} \sqcup K_{1,2},  \label{K1} \\
&&
\hspace*{.7cm}
K_{1,1} = \set{\twomatrix{\alp}{0}{0}{\alp^{*-1}}\twomatrix{1}{v/\sqrt{\eps} } {u\sqrt{\eps}}{1+uv} }{ \alp \in \calO_{k'}^\times, \; u, v \in \calO_k}, \nonumber \\
&&
\hspace*{.7cm}
K_{1,2} = \set{ \twomatrix{\alp}{0}{0}{\alp^{*-1}}\twomatrix{\pi u \sqrt{\eps}}{1+\pi uv}{1}{v/\sqrt{\eps} } }
{ \alp \in \calO_{k'}^\times, \; u, v \in \calO_k},   \nonumber
\end{eqnarray*}
and $vol(K_{1,1}) = \frac{1}{1+q^{-1}}$ and $vol(K_{1,2}) = \frac{q^{-1}}{1+q^{-1}}$ with respect to the measure on $K_1$ normalized by $vol(K_1) = 1$.

\slit
\noindent
(1) The case $x_\lam = Diag(\pi^\lam, \pi^{-\lam})$ with $\lam  \geq 0$.  

For $h = \twomatrix{1}{v/\sqrt{\eps}} {u\sqrt{\eps}}{1+uv} \in K_{1,1}$, we have
\begin{eqnarray*}
d_1(h\cdot x_\lam) &=& 
-\pi^{\lam}u^2 \eps + \pi^{-\lam}(1+uv)^2 = \pi^{-\lam} N(1+uv-\pi^\lam u\sqrt{\eps})\\
&=&
\pi^{-\lam} N(1+uv+\pi^{\lam} - 2\pi^{\lam}u \frac{1+\sqrt{\eps}}{2} ).
\end{eqnarray*}
If $u \in \frp$, then $1+uv +\pi^\lam u \in \calO_k^\times$. For $u \in \calO_k^\times$ and $r > 0$, we have
\begin{eqnarray*}
&&
vol(\set{v \in \calO_k}{1+uv + \pi^\lam u \equiv 0 (\frp^r)}) = 
vol(\set{v \in \calO_k}{ v + \pi^\lam \equiv -u^{-1} (\frp^r)}) = q^{-r},\\
\end{eqnarray*}
and for $r \geq 0$
\begin{eqnarray*}
vol(\set{(u,v) \in \calO_k^\times \times \calO_k}{v_\pi(1+uv + \pi^{\lam}u) = r }) = (1-q^{-1})^2q^{-r}.
\end{eqnarray*}
Hence we see
\begin{eqnarray}
&&
\int_{K_{1,1}} \abs{d_1(h\cdot x_\lam)}^s dh \nonumber \\
&=&
\frac{q^{\lam s}}{1+q^{-1}} \left(q^{-1} + \sum_{r=0}^{e+\lam-1} (1-q^{-1})^2q^{-r-2rs} + \sum_{r \geq e+\lam} (1-q^{-1})^2q^{-r-2(e+\lam)s} \right) \nonumber\\
&=&  \label{int-1}
\frac{q^{\lam s}}{1+q^{-1}} \left(q^{-1} + 
\frac{(1-q^{-1})^2(1-q^{-(e+\lam)-2(e+\lam)s})}{1-q^{-1-2s}} + (1-q^{-1})q^{-(e+\lam)(1+2s)} \right) .  
 \end{eqnarray}
On the other hand, for $h = \twomatrix{\pi u \sqrt{\eps}}{1+\pi uv}{1}{v/\sqrt{\eps} } \in K_{1,2}$, we have
\begin{eqnarray*}
d_1(h\cdot x_\lam) &=& 
\pi^\lam - \frac{\pi^{-\lam}v^2}{\eps}  = -\frac{\pi^{-\lam}}{\eps} \cdot N(v+\pi^\lam \sqrt{\eps})\\
&=&
-\frac{\pi^{-\lam}}{\eps} \cdot N(v-\pi^{\lam}+2\pi^\lam \frac{1+\sqrt{\eps}}{2}).
\end{eqnarray*}
Hence we see
\begin{eqnarray}
&&
\int_{K_{1,2}} \abs{d_1(h\cdot x_\lam)}^s dh \nonumber \\
&=&
\frac{q^{-1+\lam s}}{1+q^{-1}} \left( \sum_{r=0}^{e+\lam-1} (1-q^{-1})q^{-r-2rs} + \sum_{r \geq e+\lam} (1-q^{-1})q^{-r-2(e+\lam)s} \right) \nonumber\\
&=& \label{int-2}
\frac{q^{-1+\lam s}}{1+q^{-1}} \left( \frac{(1-q^{-1})(1-q^{-(e+\lam)-2(e+\lam)s})}{1-q^{-1-2s}} + q^{-(e+\lam)-2(e+\lam)s} \right).
\end{eqnarray}
By (\ref{int-1}) and (\ref{int-2}), we obtain for $s = -z - \frac12 \in \C$ 
\begin{eqnarray}
\omega(x_\lam; s) &=&
\frac{q^{-\lam z-\frac{\lam}{2}} }{ (1+q^{-1})(1-q^{2z})} \left(1-q^{2z-1}+q^{2(e+\lam)z-1}-q^{2(e+\lam+1)z} \right) \nonumber \\
&=&
\frac{ q^{-\lam z-\frac{\lam}{2}} }{ (1+q^{-1})} 
\left( \frac{q^{-\lam z}(1-q^{2z-1})}{1-q^{2z}} 
+ \frac{q^{(2e+\lam) z-1}-q^{(2e+\lam+2)z}}{1-q^{2z}} \right)\nonumber\\
&=& \label{eq:2-Case1}
\frac{q^{-\frac{\lam}{2}} q^{ez}}{ 1+q^{-1} }\left( \frac{q^{-(\lam +e)z}(1-q^{2z-1})}{1-q^{2z}} 
+ \frac{q^{(\lam + e)z}(1-q^{(-2z-1})}{1-q^{-2z}} \right).
\end{eqnarray}
 
\mslit
\noindent
(2) The case $x_\lam = \twomatrix{\pi^{\lam}(1-\eps)}{-\sqrt{\eps}}{\sqrt{\eps}}{\pi^{-\lam}}$ with $-e \leq \lam < 0$, only when $e>0$. 

Set $r = -\lam$, then $1 \leq r \leq e$.
For $h = \twomatrix{1}{v/\sqrt{\eps}} {u\sqrt{\eps}}{1+uv} \in K_{1,1}$, we have
\begin{eqnarray*}
d_1(h\cdot x_\lam) &=& 
\left( \pi^{-r} u(1-\eps)\sqrt{\eps}+(1+uv)\sqrt{\eps}\right)(-u\sqrt{\eps}) + \left(-u\eps + \pi^r(1+uv)\right)(1+uv) \nonumber\\
&=&
\pi^{-r}\left(\pi^{2r}(1+uv)^2 - 2\pi^r (1+uv)u\eps + u^2(\eps-1)\eps \right) \nonumber\\
&=&
\pi^{-r}\cdot N(\pi^r(1+uv)-u\eps -u \sqrt{\eps}) \nonumber\\
&=&
\pi^{-r}\cdot N(\pi^r(1+uv)+u(1-\eps) -2u \frac{1+\sqrt{\eps}}{2}).
\end{eqnarray*}
Since $e > 0$ and $v_\pi(1-\eps) = 2e$, we see
\begin{eqnarray*}
v_\pi(N(\pi^r(1+uv)+u(1-\eps) -2u \frac{1+\sqrt{\eps}}{2})) = 2\min\{r+v_\pi(1+uv), e+v_\pi(u)\}.
\end{eqnarray*}
We have
\begin{eqnarray*}
&&
vol(\set{(u,v) \in \calO_k^2}{v_\pi(1+uv) = 0}\\
&=&
 vol(\set{(u,v) \in \calO_k^2}{uv \in \frp}) + vol(\set{(u,v) \in \calO_k^{\times 2}}{v {\not\equiv}-u^{-1} (\frp)}) \\
&=&
2q^{-1}-q^{-2}+(1-q^{-1})(1-2q^{-1}) = 1-q^{-1}+q^{-2},
\end{eqnarray*}
and, for $j > 0$, 
\begin{eqnarray*}
&&
vol(\set{(u,v) \in \calO_k^2}{v_\pi(1+uv) = j}) = vol(\set{(u,v) \in \calO_k^{\times 2}}{v_\pi(1+uv) = j})\\
&=&
(1-q^{-1})^2q^{-j}.
\end{eqnarray*} 
Hence, for $\lam = -r < 0$, 
\begin{eqnarray}
\lefteqn{\int_{K_{1,1}} \abs{d_1(h \cdot x_\lam)}^s dh} \nonumber\\
 &=&
\frac{q^{rs}}{1+q^{-1}}\left((1-q^{-1}+q^{-2})q^{-2rs} + \sum_{j=1}^{e-r-1}(1-q^{-1})^2q^{-j-2(r+j)s} + \sum_{j \geq e-r} (1-q^{-1})^2q^{-j-2es}\right)  \nonumber\\
&=&
\frac{q^{rs}}{1+q^{-1}} \left((1-q^{-1}+q^{-2})q^{-2rs} + \frac{(1-q^{-1})^2(q^{-1-2(r+1)s}-q^{-(e-r)-2es})}{1-q^{-1-2s}} + (1-q^{-1})q^{-(e-r)-2es} \right).
 \nonumber
\\ \label{excep-1}
\end{eqnarray}
On the other hand, for $h = \twomatrix{\pi u \sqrt{\eps}}{1+\pi uv}{1}{v/\sqrt{\eps} } \in K_{1,2}$, we have
\begin{eqnarray*}
d_1(h\cdot x_\lam) &=& 
(\pi^{-r}(1-\eps)+v) - (-\sqrt{\eps}+\pi^r v/\sqrt{\eps})v/\sqrt{\eps} \\
&=&
-\frac{\pi^{-r}}{\eps} \left(\pi^{2r}v^2 - 2\pi^r v\eps + \eps(\eps-1) \right)\\
&=&
-\frac{\pi^{-r}}{\eps}\cdot N(\pi^rv-\eps-\sqrt{\eps}) = 
-\frac{\pi^{-r}}{\eps}\cdot N(\pi^rv +1-\eps -2 \frac{1+\sqrt{\eps}}{2}),\\
v_\pi(d_1(h\cdot x_\lam)) &=& 
-r+2\min\{v_\pi(v)+r, e\}.
\end{eqnarray*}
Hence we obtain 
\begin{eqnarray}
\lefteqn{\int_{K_{1,2}} \abs{d_1(h \cdot x_\lam)}^s dh} \nonumber\\
 &=&
\frac{q^{-1+rs}}{1+q^{-1}}\left(\sum_{j=0}^{e-r-1}(1-q^{-1})q^{-j-2(r+j)s} + q^{-(e-r)-2es}\right)  \nonumber\\
&=&
\frac{q^{-1+rs}}{1+q^{-1}}\left( \frac{(1-q^{-1})(q^{-2rs}-q^{-(e-r)-2es})}{1-q^{-1-2s}} + q^{-(e-r)-2es}   \right). \label{excep-2}
\end{eqnarray}

By (\ref{excep-1}) and (\ref{excep-2}), we obtain for $\lam=-r$ with  $1 \leq r \leq e$ and $s = -z-\frac12  \in \C$,  
\begin{eqnarray}
\omega(x_\lam; s) \nonumber 
&=&
\frac{q^{-rz-\frac{r}{2}} }{ (1+q^{-1})(1-q^{2z})} 
\left(q^{2rz+r} - q^{2(r+1)z+r-1} + q^{2ez+r-1} - q^{2(e+1)z+r} \right) \nonumber \\
&=&
\frac{ q^{ez+\frac{r}{2}} }{ (1+q^{-1})(1-q^{2z})} 
\left(q^{(r-e)z} - q^{(r-e+2)z-1} + q^{(-r+e)z-1} - q^{(-r+e+2)z} \right) \nonumber \\
&=&  \label{eq:2-Case2}
\frac{q^{-\frac{\lam}{2}} q^{ez} }{ 1+q^{-1}} 
\left(\frac{q^{-(\lam+e)z}(1-q^{2z-1})}{1-q^{2z}} +  \frac{q^{(\lam+e)z}(1-q^{-2z-1})}{1-q^{-2z}} \right) .
\end{eqnarray}

We have established the explicit formula of $\omega(x;s)$ by \eqref{eq:2-Case1} and \eqref{eq:2-Case2}, from which the property \eqref{eq:fun-eq m=2} follows.
\qed

\bigskip
\noindent
{\bf 2.3.} 
In this subsection we study $\omega(x; s)$ for $(m,n) = (3, 1)$ under the assumption $e \leq 1$ and show the following. The odd residual case ($e=0$) has been proved in \cite[Proposition~2.3]{HK2}. As a result, we see again the set $\xx_1^{od} = \set{x_\lam}{\lam\in\wt{\Lam_1^+}}$ forms a set of complete representatives of $K_1\backslash X_1$, since $\omega(x_\lam;z)$ takes different value for each $\lam$ for generic $z$.

\begin{prop} \label{prop: odd n=1}
Assume $e \leq 1$. Then, for $x_\lam \in \xx_1^{od}$, one has 
\begin{eqnarray} \label{eq:explicit m=3}
\omega(x_\lam; z)  
= \frac{\sqrt{-1}^\lam q^{-\lam}q^{ez}(1-q^{-1+2z})}{(1+q^{-3})(1+q^{2z})} \times \left(\frac{q^{-(\lam+e)z}(1+q^{-2+2z})}{1-q^{2z}} + 
\frac{q^{(\lam+e)z}(1+q^{-2-2z})}{1-q^{-2z}} \right). \nonumber\\
\label{eq:sph m=3}
\end{eqnarray}
In particular, for any $x \in X_1^{(od)}$, 
\begin{eqnarray} \label{eq:fun-eq m=3}
\frac{q^{-ez}(1+q^{2z})}{1-q^{-1+2z}} \cdot \omega(x;z) \in \C[q^{z}+q^{-z}],\quad
\omega(x;z) = q^{2ez}\frac{1-q^{-1+2z}}{q^{2z}-q^{-1}} \cdot \omega(x;-z).
\end{eqnarray}
\end{prop}


\begin{rem} \label{Rem:2-4}{\rm 
We expect Proposition~\ref{prop: odd n=1} holds for every $e \geq 0$. 
If the property \eqref{eq:fun-eq m=3} holds for any $e \geq 0$, then all the statements in this paper will hold for any $e \geq0$. At the moment, because of the calculation of \eqref{eq:sph m=3}, for odd $m$, we establish only for $e = 0, 1$.
} 
\end{rem}

\slit
Recall the expression of $K_1$ given in Lemma~1.6. 
We see 

\begin{center}
\begin{minipage}{13cm}
the condition ``$b, c \in \calO_{k'}$ with $N(b)+c+c^*=0$" is equivalent to \\
``$b \in \calO_{k'}, c_1 \in \calO_k$ with $N(b)+c_1 \in 2\calO_k$", where $c = \frac{N(b)+c_1}{2}+c_1\frac{1+\sqrt{\eps}}{2}$.
\end{minipage}
\end{center}

%

\begin{lem}   \label{lem: vol of b}
We normalize the Haar measures on $k'$ by $vol(\calO_k) = 1$.\\
 {\rm (1)} For $r \in \N$ and $c_1 \in \calO_k$ with $v_\pi(c_1) <r$, one has 
\begin{eqnarray*}
&&
vol(\set{b \in \calO_{k'}}{N(b) + c_1 \in \pi^r\calO_k}) = \left\{\begin{array}{ll}
0 & \mbox{if $v_\pi(c_1)$ is odd}\\
(1+q^{-1})q^{-r} & \mbox{if $v_\pi(c_1)$ is even}.
\end{array} \right. 
\end{eqnarray*}
{\rm (2)} For any $c_1 \in \calO_k$, one has 
\begin{eqnarray*}
&&
vol(\set{b \in \calO_{k'}}{v_\pi(N(b) + c_1) =  r}) \\
&&
= \left\{\begin{array}{ll}
0 & \mbox{for odd $r < v_\pi(c_1)$ and odd $v_\pi(c_1) < r$}\\[2mm]
q^{-r} & \mbox{for even $r < v_\pi(c_1)$}\\[2mm]
q^{-(r+1)} & \mbox{for odd $r = v_\pi(c_1)$}\\[2mm]
(1-q^{-1}-q^{-2})q^{-r} & \mbox{for even $r = v_\pi(c_1)$}\\[2mm]
(1-q^{-2})q^{-r} & \mbox{for even $v_\pi(c_1) < r$}.
\end{array} \right.
\end{eqnarray*}
\end{lem}

\proof
(1) Set $S(c_1, r) = \set{b \in \calO_{k'}}{N(b)+c_1=r}$. When $v_\pi(c_1)$ is odd, $S(c_1, r) = \emptyset$ and its value is $0$. When $c_1 \in \calO_k^\times$, $S(c_1, r) \subset \calO_{k'}^\times$. Since the norm map induces the surjective group homomorphism 
$$
\calO_{k'}^\times/(\pi^r) \longrightarrow \calO_k^\times/(\pi^r),
$$
which is $((1+q^{-1})q^r : 1)$-map, we see $vol(S(c_1,r)) = (1+q^{-1})q^rq^{-2r} = (1+q^{-1})q^{-r}$. When $v_\pi(c_1) = 2t > 0$, we see
\begin{eqnarray*}
 vol(S(c_1,r)) 
&= &
 vol(\set{\pi^t \xi \in \pi^t\calO_{k'}^\times}{N(\xi) + \pi^{-2t}c_1 \in \pi^{r-2t}\calO_k}) \\
& = & q^{-2t}\cdot (1+q^{-1})q^{-(r-2t)} = (1+q^{-1})q^{-r}.
\end{eqnarray*}
As for (2), the result is clear except for the case $r = v_\pi(c_1)$ is even. We see 
\begin{eqnarray*}
&&
 vol(\set{b \in \calO_{k'}}{v_\pi(N(b) + c_1) =  r})\\
&=&
vol(\pi^{r/2+1}\calO_{k'}) + vol(\set{\pi^{r/2}\xi \in \pi^{r/2} \calO_{k'}^\times}{N(\xi)-\pi^{-r}c_1 \notin \pi\calO_k })\\
& = &
q^{-r-2} + q^{-r}(1-q^{-2}-q^{-1}(1+q^{-1})) \\
&=&
(1-q^{-1}-q^{-2})q^{-r}.
\end{eqnarray*}
\qed

\begin{lem} \label{lem: vol of K11 and K12}
By the Haar measures on $k$ and $k'$ normalized by $vol(\calO_k) = vol(\calO_{k'}) = 1$, one has
\begin{eqnarray*} 
&&
vol(\set{(b,c_1) \in \calO_{k'} \times \calO_k}{N(b)+c_1 \in 2\calO_k}) = q^{-e}, \\
&&
vol(\set{(b,c_1) \in \pi\calO_{k'} \times \pi\calO_k}{N(b)+c_1 \in 2\pi\calO_k}) = q^{-(3+e)},
\end{eqnarray*}
and by the Haar measure on $G_1$ normalized by $vol(K_1) =1$, one has 
\begin{eqnarray*} 
&&
vol(K_{1,1}) = \frac{1}{1+q^{-3}}, \quad vol(K_{1,2}) = \frac{q^{-3}}{1+q^{-3}}.
\end{eqnarray*}
\end{lem}

\proof 
By Lemma~\ref{lem: vol of b}, we obtain
\begin{eqnarray*}
\lefteqn{ 
vol(\set{(b,c_1) \in \calO_{k'} \times \calO_k}{N(b)+c_1 \in \pi^e \calO_k}) \nonumber }\\
&=&
(1-q^{-1})(1+q^{-1})q^{-e} + \sum_{t =1}^{[\frac{e-1}{2}]}\, (1-q^{-1})q^{-2t} (1+q^{-1})q^{-e} + q^{-e} 
\left(\begin{array}{ll} 
q^{-e} & \mbox{if } 2\mid e\\
q^{-(e+1)} & \mbox{if } 2 \not{\mid} e
\end{array} \right) \nonumber \\
&=&
(1-q^{-2})q^{-e} + q^{-e}
\left(\begin{array}{ll} 
q^{-2} - q^{-e} & \mbox{if } 2\mid e\\
q^{-2} - q^{-(e+1)} & \mbox{if } 2 \not{\mid} e
\end{array} \right) + q^{-e}\left(\begin{array}{ll} 
q^{-e} & \mbox{if } 2\mid e\\
q^{-(e+1)} & \mbox{if } 2 \not{\mid} e
\end{array} \right) \nonumber\\
&=&
q^{-e},\\[2mm]
\lefteqn{ 
vol(\set{(b,c_1) \in \pi\calO_{k'} \times \pi\calO_k}{N(b)+c_1 \in \pi^{e+1} \calO_k}) \nonumber}\\
&=&
\sum_{t =1}^{[\frac{e}{2}]}\, (1-q^{-1})q^{-2t} (1+q^{-1})q^{-(e+1)} + q^{-(e+1)} 
\left(\begin{array}{ll} 
q^{-(e+2)} & \mbox{if } 2\mid e\\
q^{-(e+1)} & \mbox{if } 2 \not{\mid} e
\end{array} \right)\nonumber\\
&=&
q^{-(e+1)}
\left(\begin{array}{ll} 
q^{-2} - q^{-(e+2)} & \mbox{if } 2\mid e\\
q^{-2} - q^{-(e+1)} & \mbox{if } 2 \not{\mid} e
\end{array} \right) + q^{-(e+1)}\left(\begin{array}{ll} 
q^{-(e+2)} & \mbox{if } 2\mid e\\
q^{-(e+1)} & \mbox{if } 2 \not{\mid} e
\end{array} \right)\nonumber\\
&=&
q^{-(e+3)}.
\end{eqnarray*}
Now we see the volume of $K_{1,1}$ and $K_{1,2}$ as above, by the explicit description of $K_1$ in Lemma~\ref{lem: explicit K1}. 
\qed

\bigskip
In the rest of this subsection we assume $e = 1$, i.e., $2$ is a prime element in $k$. As for the calculation of $\abs{d_1(h\cdot x_\lam)}, \; \lam \geq -1$ , only the third row of $h \in K_1$ is concerned. 

\mslit
\noindent
(1) The case $\lam \geq 0$.

\medskip
For $h = \begin{pmatrix}1  &&\\ b & 1 &  \\ c& -b^* & 1 \end{pmatrix}
\begin{pmatrix}  1& d & f\\ &1&-d^*\\  &  & 1 \end{pmatrix} \in K_{1,2}$,
we have
\begin{eqnarray*}
d_1(h\cdot x_\lam) 
&=& 
\pi^\lam N(c) + N(cd-b^*) + \pi^{-\lam}N(cf-b^*d^*+1)\\
&=&
\pi^{-\lam}\left(\pi^{2\lam} N(c) + \pi^\lam N(cd-b^*) + N(1+cf-b^*d^*)\right) \in \pi^{-\lam}\calO_k^\times,
\end{eqnarray*}
hence we obtain
\begin{eqnarray} \label{eq:value of K12}
\int_{K_{1,2}} \abs{d_1(h\cdot x_\lam}^s dh = \frac{q^{-3+\lam s}}{1+q^{-3}}.
\end{eqnarray}

\bigskip
For $h = \begin{pmatrix}  &&1\\&1&-b^*\\ 1&b & c \end{pmatrix} \in K_{1,1}$, 
we have
\begin{eqnarray*}
d_1(h\cdot x_\lam) &=& 
\pi^\lam + bb^* + \pi^{-\lam}cc^* = \pi^{-\lam}(\pi^{2\lam}+\pi^\lam N(b)+N(c)). \nonumber\\
\end{eqnarray*}
Here, since we have
\begin{eqnarray}
&& \label{eq:c-c0c1}
c =  c_0+c_1\frac{1+\sqrt{\eps}}{2}  
= -\frac{N(b)}{2}+\frac{c_1}{2} \sqrt{\eps} \quad (c_0, \; c_1 \in \calO_k), \nonumber \\
&&   \label{eq:N(c)}
N(c) =
\frac{1}{4}(N(b)^2-c_1^2 \eps),
\end{eqnarray}
we see
\begin{eqnarray*}
d_1(h\cdot x_\lam) &=& 
\pi^{-\lam}\left( (\pi^\lam+\frac{N(b)}{2})^2 - \frac{c_1^2 \eps}{4} \right) \nonumber\\
&=&
\pi^{-\lam}\cdot N(\pi^\lam+\frac{N(b)}{2} - c_1 \frac{\sqrt{\eps}}{2}) \nonumber \\
&=&
\pi^{-\lam}\cdot N(\pi^\lam+\frac{N(b)+c_1}{2} - c_1 \frac{1+\sqrt{\eps}}{2} ).  \label{d1-lam}
\end{eqnarray*}
By Lemma~\ref{lem: vol of K11 and K12}, we have 
\begin{eqnarray*}
\int_{K_{1,1}} \abs{d_1(h\cdot x_\lam)}^s dh = \frac{q^{1+\lam s}}{1+q^{-3}} \sum_{r \geq 0}\, \mu(\lam,r)q^{-2rs},
\end{eqnarray*}
where
\begin{eqnarray*} \label{ylam}
&&
\mu(\lam,r) = vol(\set{(b,c_1)\in \calO_{k'}\times \calO_k}{N(b)+c_1 \in 2\calO_k, \; v_\pi(y_\lam) = r}),\\
&&
y_\lam = \pi^\lam + \frac{N(b)+c_1}{2}-c_1\frac{1+\sqrt{\eps}}{2}. \nonumber
\end{eqnarray*}
%
%
For simplicity of notation, we set $t = q^{-1}$ and $X = q^{-s}$. 
We calculate the value $\mu(\lam, r)$ case by case by using Lemma~\ref{lem: vol of b}. 
 
\noindent
the case even $r$ with $0 \leq r \leq \lam-1$: 
\begin{eqnarray*}
\mu(\lam, r)  &=&  
vol(\set{(b,c_1)}{v_\pi(c_1) = r, \; v_\pi(N(b)+c_1) \geq r+1}), \nonumber \\
&\quad +&
 vol(\set{(b,c_1)}{v_\pi(c_1) = v_\pi(N(b)+ c_1) = r+1}), \nonumber\\
&=& 
(1-t)t^r\cdot (1+t)t^{r+1} + (1-t)^{r+1}t^{r+2} =  (1-t^3)t^{2r+1}.
\end{eqnarray*}

\noindent
the case odd $r$ with $0 \leq r \leq \lam-1$: 
\begin{eqnarray*}
\mu(\lam, r)  &=&  vol(\set{(b,c_1)}{v_\pi(c_1) = v_\pi(N(b)+ c_1) = r+1}), \nonumber \\
&\quad +&
 vol(\set{(b,c_1)}{v_\pi(c_1) \geq r+2,\; v_\pi(N(b)+c_1) = r+1}), \nonumber\\
&=& 
(1-t)t^{r+1}(1-t-t^2)t^{r+1} + (1-t^2)t^{2r+3} = (1-t)t^{2r+2}
\end{eqnarray*}

\noindent
the case $r = \lam$ is even: 
\begin{eqnarray*}
\mu(\lam, \lam)  &=&  
vol(\set{(b,c_1)}{v_\pi(c_1) = \lam, \; v_\pi(N(b)+c_1) \geq \lam+1}), \nonumber \\
&\quad +&
 vol(\set{(b,c_1)}{v_\pi(c_1) = v_\pi(N(b)+ 2\pi^\lam + c_1) = \lam+1}), \nonumber\\
&\quad +&
 vol(\set{(b,c_1)}{v_\pi(c_1) \geq \lam+2, \; v_\pi(N(b)+ 2\pi^\lam) = \lam+1}), \nonumber\\
&=&
(1-t^2)t^{2\lam+1} + (1-2t)t^{\lam+1}t^{\lam+2} + t^{\lam+2}t^{\lam+2} = 
 (1-t^3)t^{2\lam+1}.
\end{eqnarray*}

\noindent
the case $r = \lam$ is odd: 
\begin{eqnarray*}
\mu(\lam, \lam)  &=&  
 vol(\set{(b,c_1)}{v_\pi(c_1) = v_\pi(N(b)+ 2\pi^\lam + c_1) = \lam+1}), \nonumber\\
&\quad +&
 vol(\set{(b,c_1)}{v_\pi(c_1) \geq \lam+2, \; v_\pi(N(b)+ 2\pi^\lam) = \lam+1}), \nonumber\\
&=&
(1-2t)t^{\lam+1}(1-t-t^2)t^{\lam+1}+t^{\lam+2}(1-t^2)t^{\lam+1} + t^{\lam+2}(1-t-t^2)t^{\lam+1} \nonumber\\
&=&
(1-t)t^{2\lam+2}
\end{eqnarray*}

\noindent
the case $r = \lam+1$ and $\lam$ is even: 
\begin{eqnarray*}
\mu(\lam, r)  &=&  
vol(\set{(b,c_1)}{v_\pi(c_1) = \lam+1, \; v_\pi(N(b)+2\pi^\lam + c_1) \geq \lam+2}), \nonumber \\
&\quad +&
 vol(\set{(b,c_1)}{v_\pi(c_1) \geq \lam+2, \; v_\pi(N(b)+ 2\pi^\lam + c_1) = \lam+2}), \nonumber\\
&=&
t^{\lam+2}t^{\lam+2}+ 0 = t^{2\lam+4}. 
\end{eqnarray*}

\noindent
the case $r = \lam+1$ and $\lam$ is odd: 
\begin{eqnarray*}
\mu(\lam, r)  &=&  
vol(\set{(b,c_1)}{v_\pi(c_1) = \lam+1, \; v_\pi(N(b)+2\pi^\lam + c_1) \geq \lam+2}), \nonumber \\
&\quad +&
 vol(\set{(b,c_1)}{v_\pi(c_1) \geq \lam+2,\; v_\pi(N(b)+2\pi^\lam+c_1) = \lam+2}), \nonumber\\
&=&
(t^{\lam+2}t^{\lam+3}+(1-2t)t^{\lam+1}(1+t)t^{\lam+2})+t^{\lam+2}(1-t^2)t^{\lam+2}\nonumber\\
&=&
(1-t^2-t^3)t^{2\lam+3}.
\end{eqnarray*}

\noindent
$\mu(\lam, r) = 0$ if $r \geq \lam+2$ and $\lam$ is even.   

\medskip
\noindent
the case $r \geq \lam+2$ and $\lam$ is odd: 
\begin{eqnarray*}
\mu(\lam, r)  &=&  
vol(\set{(b,c_1)}{v_\pi(c_1) = r, \; v_\pi(N(b)+2\pi^\lam +c_1) \geq r+1}), \nonumber \\
&\quad +&
 vol(\set{(b,c_1)}{v_\pi(c_1) \geq r+1, \; v_\pi(N(b)+2\pi^\lam + c_1) = r+1}), \nonumber\\
&=&
(1-t)t^r(1+t)t^{r+1}+t^{r+1}(1-t^2)t^{r+1}\nonumber\\
&=&
(1+t)(1-t^2)t^{2r+1}.
\end{eqnarray*}

\medskip
\noindent
By these data, we obtain the value $\int_{K_{1,1}}\abs{d_1(h\cdot x)}^sdh$ as follows:
\begin{eqnarray*}
&&
\mbox{If $\lam$ is even,} \nonumber\\[2mm]
&&
\quad \frac{X^{-\lam}}{1+t^3}\left(
   \frac{(1-t^3)(1-t^{2\lam+4}X^{2\lam+4})}{1-t^4X^4}+
   \frac{(1-t)(t^3X^2-t^{2\lam+3}X^{2\lam+2})}{1-t^4X^4}+t^{2\lam+3}X^{2\lam+2}\right); \nonumber \\
&&
\mbox{if $\lam$ is odd,} \nonumber\\[2mm]
&&
\quad \frac{X^{-\lam}}{1+t^3}\left(
   \frac{(1-t^3)(1-t^{2\lam+2}X^{2\lam+2})}{1-t^4X^4}+
   \frac{(1-t)(t^3X^2-t^{2\lam+5}X^{2\lam+4})}{1-t^4X^4}\right.\nonumber\\
&&
\qquad \qquad \left. +(1-t^2-t^3)t^{2\lam+2}X^{2\lam+2}
+ \frac{(1+t)(1-t^2)t^{2\lam+4}X^{2\lam+4})}{1-t^2X^2} \right).
\end{eqnarray*}
Together with \eqref{eq:value of K12}, we continue the calculation, where we recall the relation $s = -z-1+\frac{\pi\sqrt{-1}}{2\log q}$,\; $t = q^{-1}$ and $X = q^{-s}$. 
If $\lam$ is even, 
\begin{eqnarray} 
\omega(x_\lam;s)&=&
\frac{X^{-\lam}}{(1+t^3)(1-t^4X^4)}\left\{(1-t^3)+(1-t)t^3X^2+t^3-t^7X^4  \right. \nonumber\\
&&
\quad \left. + t^{2\lam+4}X^{2\lam+2}-(1-t^3)t^{2\lam+4}X^{2\lam+4}-c^{2\lam+7}X^{2\lam+6} \right\} \nonumber \\
&=&
\frac{(1+t^3X^2)X^{-\lam}}{(1+t^3)(1-t^4X^4)}\left\{1-t^4X^2 + c^{2\lam+4}X^{2\lam+2}(1-X^2) \right\} \nonumber\\
&=&
\frac{(1-q^{-1+2z})(\sqrt{-1})^\lam q^{-\lam-\lam z}}{(1+q^{-3})(1-q^{4z}))}\left\{1+q^{-2+2z}-q^{-2+(2\lam+2)z}(1+q^{2+2z}) \right\} \label{even=odd}\\
&=&
\frac{(\sqrt{-1})^\lam q^{-\lam}q^z(1-q^{-1+2z})}{(1+q^{-3})(1+q^{2z}))}\left\{
\frac{q^{-(\lam+1)z}(1+q^{-2+2z})}{1-q^{2z}} + 
\frac{q^{(\lam+1)z}(1+q^{-2-2z})}{1-q^{-2z}} \right\}. \nonumber\\
\label{eq:3-Case1}
\end{eqnarray}
If $\lam$ is odd, 
\begin{eqnarray} 
\omega(x_\lam;s)&=&
\frac{X^{-\lam}}{(1+t^3)(1-t^4X^4)}\left\{(1-t^3)+(1-t)t^3X^2+t^3-t^7X^4 + \right. \nonumber\\
&&
\quad \left. + t^{2\lam+4}X^{2\lam+2}+(1-t^3)t^{2\lam+4}X^{2\lam+4}-c^{2\lam+7}X^{2\lam+6} \right\} \nonumber \\
&=&
\frac{(1+t^3X^2)X^{-\lam}}{(1+t^3)(1-t^4X^4)}\left\{1-t^4X^2-c^{2\lam+4}X^{2\lam+2}(1-X^2) \right\} \nonumber\\
&=&
\frac{(1-q^{-1+2z})(\sqrt{-1})^\lam q^{-\lam-\lam z}}{(1+q^{-3})(1-q^{4z}))}\left\{1+q^{-2+2z}-q^{-2+(2\lam+2)z}(1+q^{2+2z}) \right\},\nonumber
\end{eqnarray}
which is the same with \eqref{even=odd}, and we obtain the same expression \eqref{eq:3-Case1} for odd $\lam$. Thus we have proved the formula \eqref{eq:explicit m=3} for $e = 1$ and $\lam \geq 0$.

\bigskip
\noindent
(2) We consider the remaining case for $e=1$, i.e. $\lam = -1, \: \pi = 2$, and  
$$
x_{-1} = \begin{pmatrix}\frac{1-\eps}{2} && -\sqrt{\eps}\\
&1&\\ \sqrt{\eps} && 2 \end{pmatrix}.
$$
For $h = \begin{pmatrix}1  &&\\ b & 1 &  \\ c& -b^* & 1 \end{pmatrix}
\begin{pmatrix}  1& d & f\\ &1&-d^*\\  &  & 1 \end{pmatrix} \in K_{1,2}$,
we have
\begin{eqnarray*}
d_1(h\cdot x_{-1}) 
&=& 
\frac{1-\eps}{2}N(c)+(1+b^*d^*+cf)c^*\sqrt{\eps} +N(b-c^*d^*)\\
&& \quad
+ (-\sqrt{\eps}c+2(1+b^*d^*+cf)(1+bd+c^*f^*)\\
&=&
\frac{1-\eps}{2}N(c) + N(b-c^*d^*) + (c^*-c)\sqrt{\eps}+(b^*c^*d^*-bcd)\sqrt{\eps}\\
&& \quad
+N(c)(f\sqrt{\eps}-f^*\sqrt{\eps})+2N(1+bd+c^*f^*)
\end{eqnarray*}
Since $b, c \in \pi\calO_{k'}$ satisfying $N(b)+c+c^* = 0$, we see 
$v_\pi(c-c^*) \geq 2$ and $v_\pi(h\cdot x_{-1}) = v_\pi(2N(1+bd+c^*f^*)) = 1$. Hence
\begin{eqnarray} \label{eq:value of K12-exp}
\int_{K_{1,2}} \abs{d_1(h\cdot x_\lam}^s dh = \frac{q^{-3-s}}{1+q^{-3}}.
\end{eqnarray}
For $h = \begin{pmatrix}  &&1\\&1&-b^*\\ 1&b & c \end{pmatrix} \in K_{1,1}$, 
we have
\begin{eqnarray}
d_1(h\cdot x_{-1}) &=& 
\frac{1-\eps}{2}+(c-c^*)\sqrt{\eps}+N(b)+2N(c) \nonumber.
\end{eqnarray}
Here, since we have \eqref{eq:N(c)} and 
\begin{eqnarray*}
&&
(c-c^*)\sqrt{\eps} = c_1\eps,
\end{eqnarray*}
we see
\begin{eqnarray*}
d_1(h\cdot x_{-1}) &=& 
\frac{1}{2}\left\{ (N(b)+1)^2-(c_1-1)^2\eps \right\} \nonumber\\
&=&
\frac12 N(N(b)+1-(c_1-1)\sqrt{\eps}) \nonumber \\
&=&
\frac12 N(N(b)+c_1 - 2(c_1-1)\frac{1+\sqrt{\eps}}{2}).
\end{eqnarray*}
By Lemma~\ref{lem: vol of K11 and K12}, we have 
\begin{eqnarray*}
\int_{K_{1,1}} \abs{d_1(h\cdot x_{-1})}^s dh = \frac{q}{1+q^{-3}} \sum_{r \geq 1}\, \mu(r)q^{-(2r-1)s},
\end{eqnarray*}
where
\begin{eqnarray*} \label{y}
&&
\mu(r) = vol(\set{(b,c_1)\in \calO_{k'}\times \calO_k}{N(b)+c_1 \in 2\calO_k, \; v_\pi(y) = r}),\\
&&
y = N(b)+c_1-2(c_1-1)\frac{1+\sqrt{\eps}}{2}. \nonumber
\end{eqnarray*}
Then we have
\begin{eqnarray*}
\mu(1)&=&
vol(\set{(b,c_1)}{c_1 \in \calO_k^\times, \; c_1 \notin 1+\pi\calO_k, \; v_\pi(N(b)+c_1) \geq 1}) \nonumber\\
&& \quad +
vol(\set{(b,c_1)}{c_1 \in 1+\pi\calO_k, \; v_\pi(N(b)+c_1) =1})\nonumber\\
&& \quad +
vol(\set{(b,c_1)}{c_1 \in \pi\calO_k, \; b \in \pi\calO_{k'}})\nonumber\\
&=&
(1-2t)(1+t)t + t (1-t^2)t + tt^2 \nonumber\\
&=&
(1-t^2-t^3)t;
\end{eqnarray*}
and for $r \geq 2$, 
\begin{eqnarray*}
\mu(r)&=&
vol(\set{(b,c_1)}{c_1 \in 1+\pi^{r-1}\calO_k^\times, \; v_\pi(N(b)+c_1) \geq r}) \nonumber\\
&& \quad 
+ vol(\set{(b,c_1)}{c_1 \in 1+\pi^r\calO_k, \; v_\pi(N(b)+c_1)= r}) \nonumber\\
&=&
(1-t)t^{r-1}(1+t)t^2 + t^r(1-t^2)t^r \nonumber\\
&=&
(1+t)(1-t^2)t^{2r-1}.
\end{eqnarray*}
Hence we have
\begin{eqnarray*}
\int_{K_{1,1}} \abs{d_1(h\cdot x)}^s dh &=&
\frac{1}{1+t^3}\left\{(1-t^2-t^3)X + \frac{(1+t)(1-t^2)t^2X^3}{1-t^2X^3}\right\},
\end{eqnarray*}
and together with \eqref{eq:value of K12-exp}, we obtain
\begin{eqnarray}
\omega(x_{-1};s) &=& \frac{(1-t^2)(1+t^3X^2)X}{(1+t^3)(1-t^2X^2)} \nonumber\\
&=&  \label{eq:3-Case2}
\frac{(1-q^{-2})(\sqrt{-1})^{-1}q^{1+z}(1-q^{-1+2z})}{(1+q^{-3})(1+q^{2z})},
\end{eqnarray}
which coincides with \eqref{eq:explicit m=3} for $e = 1$ and $\lam=-1$.
Thus we have established the explicit formula of $\omega(x;s)$ by \eqref{eq:3-Case1} and \eqref{eq:3-Case2}, from which the property \eqref{eq:fun-eq m=3} follows.
\qed

\bigskip
\noindent
{\bf 2.4.} 
In this subsection we give the functional equation with respect to $\tau$ for general $n$.

\begin{thm}   \label{th:f-eq tau} 
Assume $e \leq 1$ if $m$ is odd. 
For general size $n$, the spherical function satisfies the functional equation 
\begin{eqnarray*}
\omega(x;z) = q^{2ez_n} \left(\begin{array}{ll}
1 & \mbox{if } m = 2n\\[2mm]
\dfrac{1-q^{-1+2z_n}}{q^{2z_n}-q^{-1}} & \mbox{if } m = 2n+1
\end{array}\right)
\times \omega(x; \tau(z)),
\end{eqnarray*}
where $\tau(z) = (z_1, \ldots, z_{n-1}, -z_n)$.
\end{thm}

For $n = 1$ the statement has been shown in Proposition~2.2 and Proposition~2.3. 
Hereafter we assume $n \geq 2$, i.e., $m \geq 4$, and set
\begin{eqnarray*} \label{eq:tau}
w_\tau = \begin{pmatrix}
1_{n-1} &&\\
&j_r\\
&&1_{n-1}
\end{pmatrix} \in K, \quad r = m-2(n-1) \in \{2,3\}.
\end{eqnarray*}
Then the standard parabolic subgroup $P$ of $G$ attached to $\tau$ is given as follows, keeping $r$ as above,
\begin{eqnarray}
P &= &B \cup B w_\tau B \nonumber\\
&=&
\set{
\begin{pmatrix}
A&&\\
&h&\\
&&j_{n-1}A^{* -1}j_{n-1}
\end{pmatrix}
\begin{pmatrix}
1_{n-1}& \alp j_r & Bj_{n-1}\\
&1_r & -\alp^* j_{n-1}\\
&&1_{n-1}
\end{pmatrix} }
{\begin{array}{l}
A \in B_{n-1}(k')\\
h \in G_1 = U(j_r)\\
\alp \in M_{n-1,r}(k')\\
B \in M_{n-1}(k')\\
B+B^*+\alp j_r\alp^* = 0
\end{array} }, \nonumber\\
\label{eq:P for tau} 
\end{eqnarray}
where $B_{n-1}(k')$ is the Borel subgroup of $GL_{n-1}(k')$ consisting of all the upper triangular matrices. Here $d_i(x)$ is a relative $P$-invariant for $1 \leq i \leq n-1$, but $d_n( x)$ is not. We enlarge the group and the space and consider the action of $P' = P \times GL_1(k')$ on $X'=X \times V$ with $V = M_{r1}(k')$: 
\begin{eqnarray*}
(p,t)\star (x,v) = (p\cdot x, \rho(p)vt^{-1}), \quad (p,t) \in P', \; (x,v) \in X',
\end{eqnarray*}
where $\rho(p) = h \in U(j_r)$ for the decomposition of $p$ as in \eqref{eq:P for tau}. Set
\begin{eqnarray*} \label{eq:g(x,v)}
g(x,v) = \det\left[ \twomatrixplus{v^*j_r}{}{}{1_n-1} \cdot x_{(n-1+r)} \right], 
\end{eqnarray*}
where $x_{(n-1+r)}$ is the lower $(n-1+r) \times (n-1+r)$-block of $x$, and the matrix inside of $[\; ]$ is of size $n$. Though we have slightly changed the definition of $g(x,v)$ when $m=2n$, we have the following similar results as in \cite{HK} and \cite{HK2}.

\begin{lem} 
{\rm (1)} The function $g(x,v)$ is a relative $P'$-invariant on $X'$ associated by the character 
\begin{eqnarray*}
P' \ni (p, t) \longmapsto N(d_{n-1}(p))N(t)^{-1} = \psi_{n-1}(p)N(t)^{-1},
\end{eqnarray*}
and satisfies $g(x, v_0) = d_n(x)$ where $v_0 = {}^t(1,0)$ or ${}^t(1,0,0)$, according to the parity of $m$. 

{\rm (2)} For $x \in X^{op}$, there is $D_1(x) \in X_1$ satisfying 
\begin{eqnarray*}
g(x,v) = (d_{n-1}(x)D_1(x))[v]. 
\end{eqnarray*}
Here, for diagonal $x$, $D_1(x) = Diag(x_n^{-1}, x_n)$ or $Diag(x_n^{-1}, x_0, x_n)$,  according to the parity of $m$, where $x_n$ is the $n$-th diagonal entry and $x_0$ is the $(n+1)$-th diagonal entry for odd $m$.
\end{lem}

By the embedding from $K_1$ to $K = K_n$ defined by
\begin{eqnarray*}
K_1 \ni h \longmapsto \wt{h} = \begin{pmatrix}1_{n-1}&&\\&h&\\&&1_{n-1} \end{pmatrix},
\end{eqnarray*}
we see 
\begin{eqnarray}
\omega(x;s) &=&
\int_{K_1}dh \int_{K}\abs{\bfd(k\cdot x)}^s dk \nonumber \\
&=&
\int_{K_1}\int_{K}\abs{\bfd(\wt{h}k\cdot x)}^s dkdh \nonumber \\
&=&
\int_{K} \prod_{i<n}\abs{d_i(k\cdot x)}^{s_i} \int_{K_1}\abs{d_n(\wt{h}k\cdot x)}^s_n dh dk. \nonumber
\end{eqnarray}
Since we obtain, for $y \in X^{op}$  
\begin{eqnarray*}
&&
d_n(\wt{h}\cdot y) = 
g((\wt{h},1)\star (y, h^{-1}v_0)) = g(y, h^{-1}v_0) = d_{n-1}(y)D_1(y)[h^{-1}v_0]\\
&=& 
d_{n-1}(y)h^{* -1}D_1(y)[v_0] = d_{n-1}(y) d_1(j_rh^{* -1}\cdot D_1(y)) 
=
d_{n-1}(y) d_1(hj_r \cdot D_1(y)), \nonumber
\end{eqnarray*}
we have
\begin{eqnarray}
\omega(x;s) &=& 
\int_{K} \prod_{i\leq n-2}\abs{d_i(k\cdot x)}^{s_i} \cdot \abs{d_{n-1}(k\cdot x)}^{s_{n-1}+s_n} \omega^{(1)}(D_1(k\cdot x); s_n) dk. \nonumber
\end{eqnarray}
Hence we obtain, for $m =2n$, by the property \eqref{eq:fun-eq m=2}, 
\begin{eqnarray*}
\omega(x;s) &=& q^{2ez_n}\omega(x; s_1, \ldots, s_{n-2}, s_{n-1}+2s_n, -s_n-1) \nonumber \\
&=&
q^{2ez_n} \omega(x; \tau(z));
\end{eqnarray*}
and for $m = 2n+1$, by the property \eqref{eq:fun-eq m=3},
\begin{eqnarray*}
&&
\frac{1+q^{2z_n}}{1-q^{-1+2z_n}} \times \omega(x; s)\\
&=&
q^{2ez_n}\frac{1+q^{-2z_n}}{1-q^{-1-2z_n}} \times \omega(x; s_1, \ldots,s_{n-2}, s_{n-1}+2s_n+2+\frac{\pi\sqrt{-1}}{\log q}, -s_n-2-\frac{\pi\sqrt{-1}}{\log q}), 
\end{eqnarray*}
thus
\begin{eqnarray*}
\omega(x;s) &=& \frac{q^{2ez_n}(1-q^{-1+2z_n})}{q^{2z_n}-q^{-1}} \omega(x; \tau(z)),
\end{eqnarray*}
which completes the proof of Theorem~\ref{th:f-eq tau}.
\qed

\bigskip
\noindent
{\bf 2.5.} To describe the functional equation with respect to $W$, we prepare some notation. Set
\begin{eqnarray*}
&&
\Sigma = \set{\pm e_i \pm e_j, \; 2e_i}{1 \leq i, j \leq n, \; i \ne j},\quad  \Sigma^+ = \Sigma_s^+ \cup \Sigma_\ell^+,\\
&&
\Sigma_s^+ = \set{e_i + e_j,  \; e_i - e_j}{1 \leq i < j \leq n}, \quad \Sigma_\ell^+ = \set{2e_i}{1 \leq i \leq n}
\end{eqnarray*}
where $e_i$ is the $i$-th unit vector in $\Z^n, \; 1 \leq i \leq n$. We note here that 
$\Sigma$ is the set of roots of $G_n^{(ev)}$ and $\Sigma \cup \set{e_i}{1 \leq i \leq n}$ is the set of roots of $G_n^{(od)}$. 
We consider the pairing
\begin{eqnarray*}
\Z^n \times \C^n \ni (t,z) \longmapsto \pair{t}{z} = \sum_{i=1}^n t_iz_i \in \C,
\end{eqnarray*}
which satisfies
\begin{eqnarray*}
\pair{\alp}{z} = \pair{\sigma(\alp)}{\sigma(z)}, \quad (\alp \in \Sigma, \; z \in \C^n, \; \sigma \in W).
\end{eqnarray*}
%

\bigskip
\begin{thm} \label{th:feq}
Assume $e \leq 1$ if $m$ is odd. The spherical function $\omega(x;z)$ satisfies the following functional equation
\begin{eqnarray*} \label{Gamma-sigma}
\omega(x; z) = \Gamma_\sigma^{(e)}(z) \cdot \omega(x; \sigma(z)),\qquad (\sigma \in W),
\end{eqnarray*}
where 
\begin{eqnarray*} \label{sigma-factors}
&&
\Gamma_\sigma^{(e)}(z) = \dprod{\alp \in \Sigma^+(\sigma)}\, \gamma_\alp^{(e)}(z),
\quad \Sigma^+(\sigma) = \set{\alp \in \Sigma^+}{-\sigma(\alp) \in \Sigma^+},\\ 
&&
\gamma_\alp^{(e)}(z) = \left\{\begin{array}{ll}
\dfrac{1-q^{-1+\pair{\alp}{z}}}{q^{\pair{\alp}{z}} - q^{-1}} & \mbox{if } \alp \in \Sigma_s^+\\[3mm]
q^{e\pair{\alp}{z}} & \mbox{if } \alp \in \Sigma_\ell^+, \; m=2n\\ [2mm] 
\dfrac{q^{e\pair{\alp}{z}}(1-q^{-1+\pair{\alp}{z}})}{q^{\pair{\alp}{z}} - q^{-1}} & \mbox{if } \alp \in \Sigma_\ell^+ , \; m=2n+1.
\end{array}\right. \\
\end{eqnarray*}
\end{thm}

{\it Outline of a proof.}
The Weyl group $W$ is generated by $\set{\sigma_i = (i\; i+1) \in S_n}{1 \leq i \leq n-1}$ and $\tau$.  As for the Gamma factor, we have $\Gamma_{\sigma_i}^{(e)} (z) = \gamma_{e_i-e_{i+1}}^{(e)}(z)$ by Proposition~2.1, which is independent of $e$, and $\Gamma_\tau^{(e)} (z) = \gamma_{2e_n}^{(e)} (z)$ by Theorem~\ref{th:f-eq tau}. Then, by the cocycle relation of Gamma-factors, we obtain the results.
(Of course $\Gamma_\sigma^{(0)}(z)$ is the same as $\Gamma_\sigma(z)$ in \cite{HK} or \cite{HK2}, according to the parity of $m$.)
\qed

\bigskip
The following theorem can be proved in the same way as in \cite[Theorem~2.7]{HK} based on Theorem~\ref{th:feq}, where the function $G(z)$ below is the same as in \cite{HK} or \cite{HK2}, according to the parity of $m$.  

\bigskip
\begin{thm} \label{th:W-inv} 
Assume $e \leq 1$ if $m$ is odd. The function $q^{-\pair{e}{z}} G(z) \cdot \omega(x; z)$ is holomorphic
 on $\C^n$ 
and $W$-invariant, in particular it is an element in $\C[q^{\pm z_1}, \ldots, q^{\pm z_n}]^{W}$. Here $\pair{e}{z} = e(z_1+\cdots +z_n)$ and
\begin{eqnarray*}
G(z) = \prod_{\alp}\, \frac{1 + q^{\pair{\alp}{z}}}{1 - q^{-1+\pair{\alp}{z}}},
\end{eqnarray*}
where $\alp$ runs over the set $\Sigma_s^+$ for $m=2n$ and $\Sigma^+$ for $m=2n+1$.
In particular, each Gamma factor in Theorem~\ref{th:feq} is given as
\begin{eqnarray} \label{eq:gamma factor}
\Gamma_\sigma^{(e)}(z) = \frac{q^{\pair{e}{z}}}{G(z)} \cdot \frac{G(\sigma(z))}{q^{\pair{e}{\sigma(z)}}},  
\quad \sigma \in W.
\end{eqnarray}
\end{thm}

\vspace{2cm}
\Section{The explicit formula for $\omega(x;z)$}

As for the explicit formula of $\omega(x; z)$, it suffices to determine at a representative of each $K$-orbit, hence at $x_\lam, \lam \in \wt{\Lam_n^+}$  by Theorem~1.1-(1).

\begin{thm} \label{th:explicit}
Assume $ e \leq 1$ if $m$ is odd. For each $\lam \in \wt{\Lam_n^+}$, one has the explicit formula
\begin{eqnarray*} 
\omega(x_\lam; z)
 &=&
 c_n \,q^{\pair{\lam}{z_0}} \cdot \frac{q^{\pair{e}{z}}}{G(z)} \cdot Q_{\lam+e}(z; \{t\}),
\end{eqnarray*}
where $\lam+e = (\lam_1+e, \ldots, \lam_n+e) \in \Lam_n^+$, $G(z)$ is given in Theorem~2.10 (depending on the parity of $m$), 
$z_0 \in \C^n$ is the value in $z$-variable corresponding to ${\bf0} \in \C^n$ in $s$-variable, 
\begin{eqnarray} 
&&  \label{eq:z0}
z_{0,i} = \left\{\begin{array}{ll}
-(n-i+\frac12) + (n-i)\frac{\pi\sqrt{-1}}{\log q} & \mbox{if } m = 2n\\
-(n-i+1)+(n-i+\frac12)\frac{\pi\sqrt{-1}}{\log q}& \mbox{if } m = 2n+1,
\end{array}\right. \quad (1 \leq i \leq n), \nonumber \\
&& \label{eq:c-n}
c_n = \left\{\begin{array}{ll}
\dfrac{(1-q^{-2})^n}{w_m(-q^{-1})} & \textit{if } m = 2n\\[2mm]
\dfrac{(1+q^{-1})(1-q^{-2})^n}{w_m(-q^{-1})} & \textit{if } m = 2n+1,
\end{array} \right. \quad w_m(t) = \prod_{i=1}^m (1-t^i), \nonumber \\ 
&& \label{eq:Q-lam}
Q_\mu(z; \{t\}) = \sum_{\sigma \in W}\, \sigma\left(q^{-\pair{\mu}{z}} c(z; \{t\}) \right),\quad
c(z; \{t\}) = \prod_{\alp \in \Sigma^+}\, \frac{1 -t_\alp q^{\pair{\alp}{z}}}{1 - q^{\pair{\alp}{z}}}, \nonumber\\
&& \label{eq:t-alp} 
\{t\} = \{t_\alp\} \quad \mbox{with}\; \;  t_\alp = \left\{ \begin{array}{ll} 
-q^{-1} & \textit{if } \alp \in \Sigma_s^+ \\[2mm] 
q^{-1} & \textit{if } \alp \in \Sigma_\ell^+, \; m = 2n\\[2mm]
-q^{-2} & \textit{if } \alp \in \Sigma_\ell^+, \; m = 2n+1. \end{array}\right.
\end{eqnarray}
\end{thm}

\medskip
\begin{rem} \label{Rem 3-1}{\rm
We see the main part $Q_{\lam+e}(z; \{t\})$ of $\omega(x_\lam; z)$ is contained in $\calR = \C[q^{\pm z_1}, \ldots, q^{\pm z_n}]^W$ by Theorem~2.10, 
and related to Hall-Littlewood polynomial $P_{\lam+e}(z; \{t\})$ of type $C_n$ as follows (cf. \cite{Mac}, in general): 
\begin{eqnarray} \label{P-Q}
P_\mu(z; \{t\}) = \frac{1}{W_\mu(\{t\})} \cdot Q_\mu(z; \{t\}), \quad \mu \in \Lam_n^+, 
\end{eqnarray}
where $W_\mu(\{t\})$ is the Poincar{\' e} polynomial of the stabilizer $W_\mu$ of $W$ at $\mu$, and with the present choice of $t_\alp$, it is given precisely as follows 
\begin{eqnarray} 
&&
W_\mu(\{t\}) = \frac{\wt{w}_\mu(-q^{-1})}{(1+q^{-1})^{m'}}, \qquad m' = \left[\frac{m+1}{2}\right], \nonumber \\[2mm]
&& \label{eq:w-mu}
\wt{w}_\mu(t) = 
\left\{\begin{array}{ll}
w_{n_0}(t)^2 \prod_{\ell \geq 1}w_{n_\ell}(t) & \mbox{if } n = 2m\\[2mm]
w_{n_0+1}(t)w_{n_0}(t) \prod_{\ell \geq 1}w_{n_\ell}(t) & \mbox{if } n = 2m+1, \; n_0 > 0\\[2mm]
\prod_{\ell \geq 1}w_{n_\ell}(t) & \mbox{if } n = 2m+1, \; n_0 = 0,
\end{array} \right. 
\end{eqnarray}
with $n_\ell = n_\ell(\mu) = \sharp\set{i}{\mu_i = \ell}$. It is known (cf. \cite{Mac}, \cite[Proposition~B.3]{HK}) that the set 
$\set{P_\mu(z; \{t\})}{\mu \in \Lam_n^+}$ forms an orthogonal $\C$-basis for $\calR$ for each $t_\alp \in \R, \abs{t_\alp} < 1$, and $P_{\bf0}(z; \{t\}) = 1$; and we will use this property in \S 4.  
The explicit formula can be rewritten by using $P_\mu(z; \{t\})$ as 
\begin{eqnarray} \label{eq:sph by Plam}
\omega(x_\lam; z) = \frac{(1-q^{-1})^n}{w_m(-q^{-1})}\cdot  \frac{q^{\pair{e}{z}}}{G(z)} \cdot q^{\pair{\lam}{z_0}}\, \wt{w}_{\lam+e}(-q^{-1}) \cdot P_{\lam+e}(z; \{t\}), \quad (\lam \in \wt{\Lam_n^+}).
\end{eqnarray}
} 
\end{rem} 

\medskip
\begin{rem} \label{Rem 3-2} {\rm 
The influence of the residual characteristic of the base field $k$ in the explicit formula of $\omega(x_\lam; z)$ appears as shifting $\lam+e$ in $Q_{\lam}$ or $P_{\lam}$ and the factor $q^{\pair{e}{z}}$. 

Since $\omega(x_\lam; z)$ takes a different value at each $\lam \in \wt{\Lam_n^+}$ for generic $z$, we see each $x_\lam$ represents a different $K$-orbit in $X$, which completes the Cartan decomposition of $X$ (i.e. Theorem~1.1-(2)).  As we noted in Remark~2.4, if \eqref{eq:fun-eq m=3} in Proposition~\ref{prop: odd n=1} holds for $e (>0)$, one has the explicit formula for odd size $m$ for the same $e$, and the Cartan decomposition follows also. 
} 
\end{rem}

\medskip
\begin{rem} \label{rem:z0} {\rm  The vector $z_0$ in Theorem~\ref{th:explicit} can be regarded as {\it a generalization of the dual Weyl vector} as follows, and this is the reason we changed the relation between $s$ and $z$ for $m=2n$ from that in \cite{HK}(cf. Remark~\ref{rem:shift}). We remarked about this interpretation already in \cite[Remark3.3]{HK2}. 
For $v \in \Z^n$, set
\begin{eqnarray} \label{eq:height}
\{t\}^{{\rm ht}(v)} = \prod_{\beta \in \Sigma^+}\, t_\beta^{\pair{v}{\beta^\vee}/2}, \qquad\beta^\vee=\frac{2\beta}{\pair{\beta}{\beta}}. 
\end{eqnarray}
This is the generalization of the height of roots when $v \in \Sigma$ (\cite{Mac2}), while
it can be rewritten by using $z_0$ as  
\begin{eqnarray*}
\{t\}^{{\rm ht}(v)} = q^{\pair{v}{z_0}}.
\end{eqnarray*}
} 
\end{rem}

\bigskip
We prove Theorem \ref{th:explicit} in the same way as in the case $e =0$ (\cite{HK}, \cite{HK2})
 by using a general expression formula given in \cite{French} (or in \cite{JMSJ}) of spherical functions on homogeneous spaces, which is based on functional equations of finer spherical functions corresponding to $B$-orbits in $X$ and some data depending only on the group $G$. 
We have to check the assumptions there, but it has no problem since it is independent of the residual characteristic, and we omit it.

\medskip
Recall $X^{op} = \set{x \in X}{d_i(x) \ne 0, \; 1 \leq i \leq n}$ and the Borel subgroup $B$ of $G$ consisting of the upper triangular matrices in $G$. According to the $B$-orbit decomposition 
\begin{eqnarray*} \label{factor in frX-op}
&&
X^{op} = \dsqcup{u \in \calU}\, X_u, \qquad \calU = \left( \Z/2\Z \right)^n, \\
&&
X_u = \set{x \in X^{op}}{v_\pi(d_i(x)) \equiv u_1 + \cdots + u_i\pmod{2}, \; 1 \leq i \leq n},\nonumber
\end{eqnarray*}
we define finer spherical functions 
$$
\omega_u({x}; s) = \dint{K}\, \abs{\bfd(k\cdot x)}_u^s dk, \quad
\abs{\bfd(y)}_u^s = \left\{ \begin{array}{ll} 
\prod_{i=1}^n \abs{d_i(y)}^{s_i} & \mbox{if } {y} \in X_u,\\
{} & {}\\
0 & \mbox{otherwise .}
\end{array}
\right.
$$
Then, for each $\lam \in \wt{\Lam_n^+}$ and generic $z$, we have the following identity:
\begin{eqnarray} \label{eq:vector form1}
\left( \omega_u(x_\lam; s)\right)_{u\in \calU} = {c^{-1}} \sum_{\sigma \in W}\, \gamma(\sigma(z)) B^{(e)}(\sigma, z)\left(\delta_u(x_\lam; \sigma(z) ) \right)_{u \in \calU}, 
\end{eqnarray}
where 
\begin{eqnarray}
&&
c := \Sigma_{w \in W}\, [U\sigma U : U] \quad (\mbox{$U$ is the Iwahori subgroup of $K$ associated with $B$}), \nonumber \\[1.5mm]
&&
\gamma(z) := \left\{\begin{array}{ll}
\displaystyle{\prod_{\alp \in \Sigma_s^+}}\, \frac{1-q^{-2+2\pair{\alp}{z}}}{1-q^{2\pair{\alp}{z}}} \cdot \displaystyle{\prod_{\alp \in \Sigma_\ell^+}}\, \frac{1-q^{-1+\pair{\alp}{z}}}{1-q^{\pair{\alp}{z}}} & \mbox{if } m = 2n\\[2mm]
\displaystyle{\prod_{\alp \in \Sigma_s^+}}\, \frac{1-q^{-2+2\pair{\alp}{z}}}{1-q^{2\pair{\alp}{z}}} \cdot \displaystyle{\prod_{\alp \in \Sigma_\ell^+}}\, \frac{(1+q^{-2+\pair{\alp}{z}})(1-q^{-1+\pair{\alp}{z}})}{1-q^{2\pair{\alp}{z}}} & \mbox{if } m = 2n+1,\\
\end{array} \right.  \nonumber \\
&& \label{eq:delta x}
\delta_u(x_\lam; z) := \int_{U} \abs{\bfd(\nu\cdot x_\lam)}^s_ud\nu = \abs{\bfd(x_\lam)}_u^s 
= \left\{\begin{array}{ll}
 q^{\pair{\lam}{z_0}}q^{-\pair{\lam}{z}} & \mbox{if } x_\lam \in X_u\\
 0 & \mbox{otherwise},
\end{array} \right.
\end{eqnarray}
and $B^{(e)}(\sigma,z)$ is a matrix of size $2^n$ determined by the functional equation
\begin{eqnarray*}
\left( \omega_u(x_\lam; z)\right)_{u\in \calU} = B^{(e)}(\sigma, z) \left( \omega_u(x_\lam; \sigma(z))\right)_{u\in \calU}.
\end{eqnarray*}
We note here $c$ and $\gamma(z)$ are determined by the group $G=U(j_m)$ (\cite[Theorem~4.4]{Car}) and $\gamma(z)$ coincides with $c(\lam)$ there for the character $\lam(p) = (-1)^{v_\pi(p_1\cdots p_n)} \prod_{i=1}^n \abs{N(p_i)}^{z_i}$, where $p_i$ is the $n$-th diagonal entry of $p \in B$. We don't need to calculate the constant $c$ in advance, since it is determined by the property $\omega(x;s)\mid_{s=0}\, = P_{\bf0}(z) = 1$. 
We have to be careful the second equality in \eqref{eq:delta x} especially when $\lam \notin \Lam_n^+$, where we should consider the integral associated with the decomposition
$$
U = (U\cap B) U_N, \quad U_N = \set{u \in U}{u^t \in B, \; u \equiv 1_m\pmod{(\pi)}}.
$$
We explain how $B^{(e)}(\sigma, z)$ is obtained by Theorem~2.9. 
We regard a character $\chi$ of $\calU$ as a character of $k^\times$ by $\chi(a) = \chi(v_\pi(a)), \; a \in k^\times$. Then for any $\chi \in \what{\calU}$, we see 
\begin{eqnarray*}
\sum_{u \in \calU}\chi(u)\omega_u(x;z) = \omega(x; z_\chi),
\end{eqnarray*}
where $z_{\chi, i} = z_i$ or $z_i + \frac{\pi\sqrt{-1}}{\log q}$ suitably, and by Theorem~2.9 
\begin{eqnarray*}
\omega(x; z_\chi) &=& 
\Gamma^{(e)}_\sigma(z_\chi) \, \omega(x; \sigma(z_\chi)), \qquad (\sigma \in W).
\end{eqnarray*}
We may take $\sigma\chi \in \what{\calU}$ such that $\omega(x; \sigma(z_\chi)) = \omega(x; \sigma(z)_{\sigma\chi})$, where $(\sigma\chi)(u) = \chi(\sigma^{-1}(u)), \; u \in \calU$. Thus we obtain
\begin{eqnarray}  \label{eq:vector form2}
\Big( \chi(u)\Big)_{\chi, u} \Big( \omega_u(x; z) \Big)_{u \in \calU}  = 
\wt{\Gamma^{(e)}}(\sigma, z) \Big( (\sigma\chi)(u)\Big)_{\chi, u} \Big( \omega_u(x; \sigma(z)) \Big)_{u \in \calU}, 
\end{eqnarray}
where $\wt{\Gamma^{(e)}}(\sigma, z)$ is the diagonal matrix with $\Gamma_\sigma^{(e)}(z_\chi)$ as the $\chi$-diagonal entry, and  
\begin{eqnarray} \label{eq:matrix}
B^{(e)}(\sigma, z) =  {\Big( \chi(u)\Big)_{\chi, u}}^{-1} \; \wt{\Gamma^{(e)}}(\sigma, z) \, \Big( (\sigma\chi)(u)\Big)_{\chi, u}.
\end{eqnarray}
We set the first row for $(\chi,u) \in \what{\calU}\times \calU$ as the trivial character $\bf1$. Then the first entry in the left hand side of \eqref{eq:vector form2} is equal to $\omega(x;z)$, and 
we note $z_{\bf1} =z$ and $(\sigma{\bf1})(u) = 1, \; u \in \calU$.
Hence we have by \eqref{eq:vector form1}, \eqref{eq:vector form2}, and \eqref{eq:matrix},
\begin{eqnarray}
\omega(x; z) &=& 
{c^{-1}} \, \sum_{\sigma \in W}\, \gamma(\sigma(z)) \Gamma_\sigma^{(e)}(z) q^{\pair{\lam}{z_0}}q^{-\pair{\lam}{\sigma(z)}}   \nonumber\\
&=&
\frac{q^{\pair{\lam}{z_0}}}{c} \cdot \frac{ q^{\pair{e}{z}}}{G(z)} \, \sum_{\sigma \in W}\, \sigma\left( \gamma(z) G(z) q^{-\pair{\lam+e}{z}}\right)  \quad (\mbox{by } \eqref{eq:gamma factor}) \nonumber\\
&=& \label{eq:sph-c}
\frac{q^{\pair{\lam}{z_0}}}{c} \cdot \frac{q^{\pair{e}{z}}}{G(z)} \, \sum_{\sigma \in W}\, \sigma\left(q^{-\pair{\lam+e}{z}}\, \prod_{\alp \in \Sigma^+}\, \frac{1-t_\alp q^{\pair{\alp}{z}}}{1-q^{\pair{\alp}{z}}} \right),
\end{eqnarray}
where $t_\alp$ is given as in \eqref{eq:t-alp}.
By \eqref{eq:sph-c}, we have, taking $\lam = (-e)$ and $z = z_0$, 
\begin{eqnarray*}
c &=&
\frac{1}{G(z_0)} \sum_{\sigma\in W}\, \sigma\left(\prod_{\alp}\, \frac{1-t_\alp q^{\pair{\alp}{z_0}}}{1-q^{\pair{\alp}{z_0}}} \right),
\end{eqnarray*}
which is the same as in case $e = 0$, and $c^{-1} = c_n$ in \eqref{eq:c-n}.  
\qed

\vspace{2cm}
\Section{The structure of the Schwartz space}
We keep the assumption that $e \leq 1$ if $m$ is odd.
We define the Schwartz space
\begin{eqnarray*} \label{eq:SKX}
\SKX = \set{\vphi: X \longrightarrow \C}{\mbox{left $K$-invariant, compactly supported}},
\end{eqnarray*}
and study $\hec$-module structure and Plancherel formula about it. Based on the explicit formula in \S 3, we modify the spherical function by using the value at $x_{(-e)}$ as 
\begin{eqnarray}
&& 
\Psi(x;z) =\omega(x;z)/\omega(x_{(-e)}; z) \in \calR = \C[q^{\pm z_1}, \ldots, q^{\pm z_n}]^W.
\end{eqnarray}
Then, we have
\begin{eqnarray}
&&   \label{eq:Psi-lam}
f * \Psi(x;z) = \lam_z(f) \Psi(x;z), \quad (f \in \hec),\\[2mm]
&&
\qquad \quad \lam_z: \hec \stackrel{\sim}{\longrightarrow} \calR_0 = \C[q^{\pm 2z_1}, \ldots, q^{\pm 2z_n}]^W,  \quad ({\rm cf.} \eqref{eq:Satake trans z}), \nonumber\\[2mm]
&&\label{eq:modified sph}
\Psi(x_\lam;z) = q^{\pair{\lam+e}{z_0}}\frac{\wt{w}_{\lam+e}(-q^{-1})}{\wt{w}_{\bf0}(-q^{-1})} \cdot P_{\lam+e}(z; \{t\}), \quad (\lam \in \wt{\Lam_n^+}).
\end{eqnarray} 
Here the value $\Psi(x_\lam;z)$ for dyadic case, i.e., the case $e = v_\pi(2) > 0$,  coincides with the value $\Psi(x_{\lam+e};z)$ for odd residual case. Hence all the results of this section are parallel to odd residual case (\cite[\S 4]{HK} or \cite[\S 4]{HK2}). 

\bigskip
We define the spherical Fourier transform on $\SKX$ by
\begin{eqnarray}
\begin{array}{lcll}
F: &\SKX &\longrightarrow&  \calR\\
{} & \vphi & \longmapsto & F(\vphi)(z) = \int_{X}\vphi(x)\Psi(x;z) dx,
\end{array}
\end{eqnarray}
where $dx$ is a $G$-invariant measure on $X$, and we fix the normalization of $dx$ later. 
The Hecke algebra $\hec$ acts on $\SKX$ by convolution product
\begin{eqnarray*}
f*\vphi(x) = \int_{G}f(g)\vphi(g^{-1}\cdot x) dg, \quad (f \in \hec, \; \vphi \in\SKX),
\end{eqnarray*}
where $dg$ is the Haar measure on $G$, and on $\calR$ through Satake isomorphism $\lam_z$. Then the map $F$ is compatible with $\hec$-action as follows  
\begin{eqnarray}
F(f*\vphi)(z) = \lam_z(f)\vphi(z), \qquad (f \in \hec, \; \vphi \in \SKX).
\end{eqnarray}
The space $\SKX$ is spanned by $\set{ch_\lam}{\lam \in \wt{\Lam_n^+}}$, where $ch_\lam$ is the characteristic function of $K\cdot x_\lam$, and we have by \eqref{eq:modified sph}
\begin{eqnarray} \label{eq:F-ch lam}
F(ch_\lam)(z) = q^{\pair{\lam+e}{z_0}}\frac{\wt{w}_{\lam+e}(-q^{-1})}{\wt{w}_{\bf0}(-q^{-1})} \cdot v(K\cdot x_\lam) \cdot P_{\lam+e}(z; \{t\}),
\end{eqnarray}
where $v(K\cdot x_\lam)$ the volume of $K\cdot x_\lam$ with respect to $dx$.
Since the set $\set{P_\mu(z; \{t\})}{\mu \in \Lam_n^+}$ forms a $\C$-basis for $\calR$, the map $F$ is an $\hec$-module isomorphism. Thus we see the following.

\begin{thm} \label{thm:hec-iso}
Assume $e \leq 1$ if $m$ is odd. The spherical Fourier transform $F$ gives an $\hec$-module isomorphism
$$
\SKX \stackrel{\sim}{\longrightarrow} \C[q^{\pm z_1}, \ldots, q^{\pm z_n}]^W (= \calR),
$$
where $\calR$ is regarded as $\hec$-module via $\lam_z$. In particular $\SKX$ is a free $\hec$-module of rank $2^n$.
\end{thm}

Each spherical functions on $X$ is associated with some $\lam_z$ like as \eqref{eq:Psi-lam}, and it is determined by the class of $z$ in $\left(\C\big{/}\frac{2\pi\sqrt{-1}}{\log q}\Z\right)^n \big{/}W$. The dimension of spherical functions associated with the same $\lam_z$ is at most $2^n$ by Theorem~\ref{thm:hec-iso}, and we can give a basis as below.     

\begin{cor} \label{cor:basis of sph}
Assume $e \leq 1$ if $m$ is odd. All the spherical functions on $X$ are parametrized by eigenvalues $z \in \left(\C\big{/}\frac{2\pi\sqrt{-1}}{\log q}\Z\right)^n \big{/}W$ through Satake isomorphism $\lam_z$. The set \\$\set{\Psi(x; z + u)}{u \in \{0, \frac{\pi\sqrt{-1}}{\log q} \}^n}$ forms a basis of spherical functions on $X$ corresponding to $z$.
\end{cor}

\bigskip
We will give the Plancherel formula on $\SKX$. Recall the notation $c(z;\{t\}), \, P_\mu(z;\{t\})$ and $\wt{w}_\mu(-q^{-1})$ given in Theorem~\ref{th:explicit} and Remark~\ref{Rem 3-1}. We define an inner product on $\calR$ by 
\begin{eqnarray}
\pair{P}{Q}_\calR = \int_{\fra^*} P(z)\ol{Q(z)}d\mu(z), \quad (P, Q \in \calR),
\end{eqnarray}
where
\begin{eqnarray} \label{eq:fra and meas}
&&
\fra^* = \left\{\sqrt{-1}\left(\R\Big{/}\frac{2\pi}{\log q}\Z \right) \right\}^n, \nonumber\\[2mm]
&&
d\mu(z) = \frac{1}{n! 2^n}\cdot \frac{\wt{w}_{\bf0}(-q^{-1})}{(1+q^{-1})^{m'}}\cdot  \frac{1}{\abs{c(z;\{t\})}^2}\, dz, \quad m ' = \left[\frac{m+1}{2}\right], 
\end{eqnarray}
and $dz$ is the Haar measure on $\fra^*$.  
In the following, for simplicity we write
\begin{eqnarray*}
P_{\lam+e} = P_{\lam+e}(z; \{t\}) \in \calR, \quad \wt{w}_\lam = \wt{w}_\lam(-q^{-1}) \in \R, \qquad (\lam \in \wt{\Lam_n^+})
\end{eqnarray*}
Then, by \cite[Proposition B.3]{HK}, we have 
\begin{eqnarray} \label{eq:P-lam P-mu}
\pair{P_{\lam+e}}{P_{\mu+e}}_\calR = \delta_{\lam, \mu} 
\frac{ \wt{w}_{\bf0}  } { \wt{w}_{\lam+e}}, 
\qquad (\lam, \mu \in \wt{\Lam_n^+}),
\end{eqnarray}
and by using \eqref{eq:F-ch lam}, 
\begin{eqnarray} 
\pair{F(ch_\lam)}{F(ch_\lam)}_\calR 
& = & \label{eq:value F-side}
\delta_{\lam,\mu} \, q^{2\pair{\lam+e}{\real(z_0)}}\frac{\wt{w}_{\lam+e}}{\wt{w}_{\bf0}} \cdot v(K\cdot x_\lam)^2.
\end{eqnarray}
Since there are precisely two $G$-orbits in $X$ represented by $x_0$ and $x_1$ (Theorem~1.1-(3)), we may normalize the $G$-invariant measure $dx$ on each orbit by fixing the volume of $K\cdot x_0$ and $K\cdot x_1$, where $x_0 = 1_m$ and $x_1=x_{\langle 1\rangle}$ with $\langle 1 \rangle = (1, 0, \ldots, 0)$.

\begin{lem} \label{lem:v(Kx)}
By the normalization of the $G$-invariant measure $dx$ on $X$ given as 
\begin{eqnarray} \label{eq:normalize}
v(K\cdot x_0) = q^{-2\pair{e}{\real(z_0)}}\, \frac{\wt{w}_{\bf0}(-q^{-1})}{\wt{w}_e(-q^{-1})}, \quad 
v(K\cdot x_1) = q^{-2\pair{\langle 1 \rangle+e}{\real(z_0)}}\, \frac{\wt{w}_{\bf0}(-q^{-1})}{\wt{w}_{\langle 1 \rangle + e}(-q^{-1})}, 
\end{eqnarray}
one has
\begin{eqnarray*}
v(K\cdot x_\lam) = q^{-2\pair{\lam+e}{\real(z_0)}}\, \frac{\wt{w}_{\bf0}(-q^{-1})}{\wt{w}_{\lam+e}(-q^{-1})}, \quad 
 \lam \in \wt{\Lam_n^+}.
\end{eqnarray*}
\end{lem}

\proof
For any $f \in \hec$ and $\mu \in \wt{\Lam_n^+}$, we may write
\begin{eqnarray} \label{eq:first}
f*ch_\mu = \sum_{\nu \in\wt{\Lam_n^+}} a^\mu_\nu(f)\, ch_\nu, \quad (a^\mu_\nu(f) \in \C).
\end{eqnarray}
For $f \in \hec$ and $\lam \in \wt{\Lam_n^+}$, we have 
\begin{eqnarray}
(f*\Psi(\; ;z))(x_\lam) &=& \sum_{\mu \in \wt{\Lam_n^+}} \Psi(x_\mu;z)(f*ch_\mu)(x_\lam) 
=
\sum_{\mu \in \wt{\Lam_n^+}} \Psi(x_\mu;z)a^\mu_\lam(f) \nonumber\\
&=&
\sum_{\mu \in \wt{\Lam_n^+}} q^{\pair{\mu+e}{z_0}}\frac{\wt{w}_{\mu+e}}{\wt{w}_{\bf0}}\, a^\mu_\lam(f) \cdot P_{\mu+e},  \label{LHD}
\end{eqnarray}
where we used \eqref{eq:modified sph} and the summation over $\wt{\Lam_n^+}$ is essentially a finite sum, since the support of $f$ is compact. 
On the other hand, by \eqref{eq:Psi-lam}, we have
\begin{eqnarray}
(f*\Psi(\; ;z))(x_\lam) &=& \lam_z(f) \Psi(x_\lam;z) \nonumber\\
&=&
q^{\pair{\lam+e}{z_0}}\frac{\wt{w}_{\lam+e}}{\wt{w}_{\bf0}} \cdot \lam_z(f) P_{\lam+e}.\label{RHD}
\end{eqnarray}
Taking the inner product of \eqref{LHD} and \eqref{RHD} with $P_{\mu+e}$, we have by \eqref{eq:P-lam P-mu}
\begin{eqnarray} \label{eq:M1}
q^{\pair{\mu+e}{z_0}}a^\mu_\lam(f) =
q^{\pair{\lam+e}{z_0}}\frac{\wt{w}_{\lam+e}}{\wt{w}_{\bf0}}\pair{\lam_z(f)P_{\lam+e}}{P_{\mu+e}}_\calR, 
\quad (f \in \hec, \; \lam, \mu \in \wt{\Lam_n^+}).
\end{eqnarray}
Applying the spherical transform $F$ to the both side of \eqref{eq:first}, we have 
\begin{eqnarray*} \label{eq:M2}
v(K\cdot x_\mu) q^{\pair{\mu+e}{z_0}} \frac{\wt{w}_{\mu+e}}{\wt{w}_{\bf0}}\lam_z(f) P_{\mu+e}
= \sum_{\nu\in\wt{\Lam_n^+}} a^\mu_\nu(f) v(K\cdot x_\nu)q^{\pair{\nu+e}{z_0}}\frac{\wt{w}_{\nu+e}}{\wt{w}_{\bf0}}P_{\nu+e},
\end{eqnarray*}
taking the inner product of the both side of the above identity with $P_{\lam+e}$, we obtain 
\begin{eqnarray}
v(K\cdot x_\mu) q^{\pair{\mu+e}{z_0}} \frac{\wt{w}_{\mu+e}}{\wt{w}_{\bf0}}\pair{\lam_z(f)P_{\mu+e}}{P_{\lam+e}}_\calR
=
a^\mu_\lam(f) v(K\cdot x_\lam)q^{\pair{\lam+e}{z_0}}, \nonumber &&\\
(f \in \hec, \; \lam, \mu \in \wt{\Lam_n^+}). &&\label{eq:M3}
\end{eqnarray}
Now assume $\abs{\lam}\equiv \abs{\mu} \pmod{2}$ and take $f_1 \in \hec$ to be the characteristic function of $Kg_1K$ such that $x_\lam=g_1\cdot x_\mu$. Then 
\begin{eqnarray}
a^\mu_\lam(f_1) \ne 0, \quad 
\pair{\lam_z(f_1)P_{\lam+e}}{P_{\mu+e}}_\calR = \pair{\lam_z(f_1)P_{\mu+e}}{P_{\lam+e}}_\calR,
\end{eqnarray}
Hence we obtain by \eqref{eq:M1} and \eqref{eq:M3}
\begin{eqnarray} \label{eq:M4}
\frac{v(K\cdot x_\lam)}{v(K\cdot x_\mu)} = q^{2\pair{\mu-\lam}{z_0}}\, \frac{\wt{w}_{\mu+e}}{\wt{w}_{\lam+e}} = q^{2\pair{\mu-\lam}{\real(z_0)}}\, \frac{\wt{w}_{\mu+e}}{\wt{w}_{\lam+e}}, \quad \mbox{if}\; \abs{\lam}\equiv \abs{\mu}\pmod{2}.
\end{eqnarray}
Since $x_\lam \in G\cdot x_0$ if and only if $\abs{\lam}\equiv 0\pmod{2}$ for $\lam \in \wt{\Lam_n^+}$, under the normalization of $dx$ as in \eqref{eq:normalize}, we obtain the volume $v(K\cdot x_\lam)$ by \eqref{eq:M4}, which completes the proof.  
\qed

\bigskip
We take the normalization as in Lemma~\ref{lem:v(Kx)}. Then by \eqref{eq:value F-side}, we see 
\begin{eqnarray*}
\int_{X} ch_\lam(x) \ol{ch_\mu(x)}dx
= \delta_{\lam,\mu} v(K\cdot x_\lam) =
\int_{\fra^*}F(ch_\lam)(z)\ol{F(ch_\mu)(z)} d\mu(z), \quad (\lam, \mu \in \wt{\lam_n^+}).
\end{eqnarray*}
Since $\SKX$ is spanned by the set $\set{ch_\lam}{\lam \in \wt{\Lam_n^+}}$, we obtain following theorem.

\begin{thm} {\rm (Plancherel formula on $\SKX$)} \label{th:Plancherel}
Assume $e \leq 1$ if $m$ is odd. For any $\vphi, \psi \in \SKX$, one has
\begin{eqnarray*}
\int_X\, \vphi(x)\ol{\psi(x)} dx = \int_{\fra^*}\, F(\vphi)(z) \ol{F(\psi)(z)} d\mu(z),
\end{eqnarray*}
where $dx$ is normalized as in Lemma~\ref{lem:v(Kx)}, and $\fra^*$ and $d\mu(z)$ are given in \eqref{eq:fra and meas}.
\end{thm}

\begin{cor} {\rm (Inversion formula)}\; 
Assume $e \leq 1$ if $m$ is odd. For any $\vphi \in \SKX$ and $x \in X$, one has
\begin{eqnarray*}
\vphi(x) = \int_{\fra^*}\, F(\vphi)(z) \Psi(x;z)d\mu(z).
\end{eqnarray*}
\end{cor}

\proof
For any $\vphi \in \SKX$ and $x \in X$, we have by Theorem~\ref{th:Plancherel}
\begin{eqnarray*}
\vphi(x) &=& \frac{1}{v(K\cdot x)} \int_{X} \vphi(y)\ol{ch_{K\cdot x}(y)} dx\\
&=&
\frac{1}{v(K\cdot x)} \int_{\fra^*} F(\vphi)(z) \ol{F(ch_{K\cdot x})(z)} d\mu(z)\\
&=&
\int_{\fra^*} F(\vphi)(z) \ol{\Psi(x;z)} d\mu(z)\\
&=&
\int_{\fra^*} F(\vphi)(z) \Psi(x;z) d\mu(z).
\end{eqnarray*}
\qed

\vspace{2cm}
\bibliographystyle{amsalpha}

\vspace{2cm}

\begin{flushright}
Yumiko HIRONAKA\\

\vspace{2mm}
Department of Mathematics\\ 
Faculty of Education and Integrated Sciences\\
Waseda University\\
Nishi-Waseda, Tokyo 169-8050, JAPAN

\makeatletter

\vspace{2mm}
e-mail : hironaka@waseda.jp
\end{flushright}

\end{document}